\documentclass{article}
\usepackage{amsfonts,amsmath,amssymb,amsbsy}
\usepackage{graphicx}
\usepackage{epsfig}
\usepackage{theorem}
\newtheorem{theorem}{Theorem}
\newtheorem{lemma}{Lemma}

\newtheorem{defin}{Definition}

\begin{document}

\begin{titlepage}

\title{Fractal analysis of Neimark-Sacker bifurcation}

\author{Lana Horvat Dmitrovi\' c}

\maketitle

\begin{abstract}

In this paper we show how a change of box dimension of the orbits of two-dimensional discrete dynamical systems is connected to their bifurcations in a nonhyperbolic fixed point. This connection is already shown in the case of one-dimensional discrete dynamical systems (see \cite{laho1},\cite{neveda}). Namely, at the bifurcation point the box dimension changes from zero to a certain positive value which is connected with the type of bifurcation.  First, we study a two-dimensional discrete dynamical system with only one multiplier on the unit circle, and get the result for the box dimension of the orbit on the center manifold. Then we consider the planar discrete system undergoing a Neimark-Sacker bifurcation. It is shown that box dimension depends on the order of the nondegeneracy at the nonhyperbolic fixed point and on the angle-displacement map. As it was expected, we prove that the box dimension is different in rational and irrational case.
\end{abstract}

 

\end{titlepage} 

\def\b{\beta}
\def\g{\gamma}
\def\d{\delta}
\def\l{\lambda}
\def\o{\omega}
\def\ty{\infty}
\def\e{\varepsilon}
\def\f{\varphi}
\def\:{{\penalty10000\hbox{\kern1mm\rm:\kern1mm}\penalty10000}}
\def\st{\subset}
\def\stq{\subseteq}
\def\q{\quad}
\def\M{{\cal M}}
\def\cal{\mathcal}
\def\eR{\mathbb{R}}
\def\eN{\mathbb{N}}
\def\Ze{\mathbb{Z}}
\def\Qu{\mathbb{Q}}
\def\Ce{\mathbb{C}}
\def\ov#1{\overline{#1}}
\def\D{\Delta}
\def\O{\Omega}

\def\bg{\begin}
\def\eq{equation}
\def\bgeq{\bg{\eq}}
\def\endeq{\end{\eq}}
\def\bgeqnn{\bg{eqnarray*}}
\def\endeqnn{\end{eqnarray*}}
\def\bgeqn{\bg{eqnarray}}
\def\endeqn{\end{eqnarray}}

\newcount\remarkbroj \remarkbroj=0
\def\remark{\advance\remarkbroj by1 \smallskip{Remark\ \the\remarkbroj.}\enspace\ignorespaces\,}

\section{Introduction}

In recent years, it was shown that fractal analysis can be applied to the solutions of differential equations and dynamical systems (see \cite{lapo}, \cite{pa}, \cite{pazuzu}, \cite{zuzu4}). In the case of dynamical systems, fractal analysis consists of studying the box dimension and Minkowski content of trajectories or orbits. Several articles with fractal analysis of bifurcations of dynamical systems (see \cite{neveda}, \cite{belg}, \cite{zuzu}, \cite{zuzu3}) showed that there is a direct connection between the change in box dimension of trajectories of dynamical systems and the bifurcation of that system. Around the hyperbolic singularities the box dimension is trivial ($0$), while around the nonhyperbolic singularities the box dimension is positive and connected to the appropriate bifurcation. In fact, box dimension shows the multiplicity of singularity or limit cycle in the case of weak focus or multiple limit cycle (see \cite{zuzu}). Beside this, in the article \cite{marezu}, multiplicity was also connected to the growth of $\varepsilon$- neighbourhood of the orbit near the homoclinic loop. At the beginning, this phenomenon was studied only in the continuous dynamical systems, but the same occur in the discrete systems. The article \cite{neveda} showed that the above result for box dimension can also be proven for generic saddle-node and period-doubling bifurcations of one-dimensional discrete dynamical systems. In the paper \cite{laho1} the result from the previous paper were generalized to the class of finitely nondegenerate maps in $\eR$ and applied it to one and two-parameter bifurcations with generalized sufficient conditions. So in the one-dimensional systems, the box dimension at the bifurcation point shows the maximum number of possible fixed point that can bifurcate in a given family of systems. Of course, now we are interested in this connection for other bifurcations.

Hence, it is natural to continue with the fractal analysis of two-dimensional discrete dynamical systems. For the result about the trivial box dimension of the orbit near the hyperbolic fixed point in $\mathbb{R}^{n}$, see the article \cite{laho2}. In this article we will show the results about the box dimension of the orbit around a nonhyperbolic fixed point in $\mathbb{R}^2$. We look at the following bifurcations:  one multiplier on the unit circle and two complex conjugated multipliers with the nonresonant condition. For the bifurcation of nonhyperbolic fixed point with only one multiplier on the unit circle, we get similar result for the orbit on the center manifold, which correspond to the bifurcation of appropriate one-dimensional system. Regarding the Neimark-Sacker bifurcation (see \cite{ne}, \cite{sa}, \cite{kuz}), the situation is more complicated. As in the bifurcation analysis of this bifurcation, the main problem is the difference of the dynamics between the irrational and rational case. We will see that the box dimensions for this cases are also different. Since Neimark-Sacker bifurcation is, in fact, Hopf bifurcation for maps, we will see how this result is related to the result for Hopf bifurcation showed in \cite{zuzu}. 

The main motivation for studying the box dimension of dynamical systems is related to the problem of multiplicity or cyclicity of fixed point, singularity or limit cycle. We are interested to explore how box dimension behaves during the bifurcations in $\mathbb{R}^2$ of discrete and continuous systems, its connection to the cyclicity problem and how it can be found helpful in bifurcation analysis of some complicated bifurcations in higher dimension. Moreover, this connection could be applied as numerical tool in bifurcation analysis since there exists effective algorithms for calculating box dimension. 

Now, we will introduce the notions of box dimension and Minkowski content, as basic terms which are used in our fractal analysis. First, we recall that the fractal dimension such as box dimension (also known as Minkowski dimension, Minkowski-Bouligand dimension, capacity dimension, limit capacity), can be used to analyze various objects such as various sets, graphs of a function, attractors, trajectories, etc. Fractal analysis consists of getting the  results of the fractal dimension and putting them in the context according to other properties of the studied object. In the case of dynamical systems, we connect the box dimension with the bifurcation, its type and the number of bifurcating objects. In our analysis we use the box dimension because, for the orbit of discrete dynamical system, the Hausdorff dimension fails to show anything. Namely, because of its property of countable stability, the Hausdorff dimension does not 'see' the countable sets at all. On the other hand, the box dimension is only finitely stable so it can 'see' the countable sets and the difference between the orbits around the hyperbolic and nonhyperbolic fixed point. 

Now we recall the notions of box dimension and Minkowski content. For further details see e.g. \cite{fa}, \cite{t}, \cite{zu2}.
Let $A\st\eR^N$ be bounded. The $\e$-neighbourhood of $A$ is defined by
$A_\e=\{y\in\eR^N\:d(y,A)<\e\}$.\\ Let $s\ge0$.
\textit{The lower and upper $s$-dimensional Minkowski contents of $A$} are defined by
$$\M_*^s(A):=\liminf_{\e\to0}\frac{|A_\e|}{\e^{N-s}},\,\,\,\,\M^{*s}(A):=\limsup_{\e\to0}\frac{|A_\e|}{\e^{N-s}}.$$
Then \textit{the lower and upper box dimension} are defined by 
$$\underline\dim_BA=\inf\{s>0:\M_*^s(A)=0\}, \,\,\,\,\ov\dim_BA=\inf\{s>0:\M^{*s}(A)=0\}.$$
If $\underline\dim_BA=\ov\dim_BA$ we denote it by $\dim_BA$.
If there exists $d\ge0$ such that\ $0<\M_*^d(A)\le \M^{*d}(A)<\ty,$ then we say that set $A$ is \textit{Minkowski nondegenerate}. Clearly, then $d=\dim_B A$. If $|A_\e|\simeq \e^{s}$ for $\e$ small, then $A$ is Minkowski nondegenerate set and 
$\dim_B A=N-s$. If $\M_*^s(A)=\M^{*s}(A)=\M^d(A)\in(0,\ty)$ for some $d\ge0$, then 
$A$ is said to be {\it Minkowski measurable}. Clearly, then $d=\dim_BA$.\\

Let $A\subset \mathbb{R}^{n}$ be a disjoint bounded set and $F:\mathbb{R}^{n}\rightarrow \mathbb{R}^{n}$ is a Lipschitz map. Then it holds
$$\dim_{B}F(A)\leq \dim_{B}A.$$ 
We say that $F:\Omega\rightarrow\Omega'$, where $\Omega,\Omega'$ are open sets, is a bilipschitz map if there exist positive constants $A$ and $B$ such that
$$A\left\|x-y\right\|\leq\left\|F(x)-F(y)\right\|\leq B \left\|x-y\right\|,$$
for every $x,y\in\Omega$.
If $F$ is a bilipschitz mapping, than $$\dim_{B}A=\dim_{B}F(A).$$

Box dimension has also the property of finite stability: Let $A_{i}$, $i=1,\ldots,k$ be disjointed bounded sets. Then 
\begin{equation} \label{konstab}
\dim_{B}(\bigcup_{i=1}^{k} A_{i})=\max_{i\in\{1,\ldots,k\}} \dim_{B} A_{i}.
\end{equation}
If $\dim_{B}A_{i}=d$ for every $i$, then $\dim_{B}\bigcup A_{i}=d$.

In the paper the following definitions are also used.
We say that any two sequences $(a_n)_{n\ge1}$ and $(b_n)_{n\ge1}$ of positive real numbers are {\it comparable} and write $a_n\simeq b_n$ as $n\to\ty$ if $A\le a_n/b_n\le B$ for some $A,B>0$  and  $n$ sufficiently big.
Analogously, two positive functions $f,g:(0,r)\rightarrow \eR$ are comparable and we write $f(x)\simeq g(x)$ as 
$x \rightarrow 0$ if $f(x)/g(x)\in [A,B]$ for $x$ small enough.

We will study a two-dimensional discrete dynamical system 
$$\mathbf{x}_{n+1}=\mathbf{F}(\mathbf{x}_{n}),\, \mathbf{x}_1\in \mathbf{R}^{2}$$
generated by a $C^{k}$ function $\mathbf{F}:\eR^2\mapsto \eR^2$. The \textit{orbit} of a system is a sequence $(\mathbf{x}_{n})_{n\geq1}$ such that $\mathbf{x}_{n+1}=\mathbf{F}(\mathbf{x}_{n})$ for some $\mathbf{x}_1\in \eR^2$. 
Let $x_0=0$ be a fixed point ($\mathbf{F}(\mathbf{x}_0)=\mathbf{x}_0$) of that system and let $A$ be a Jacobi matrix $DF(x_0)$ at $x_0$. The eigenvalues $\lambda_1,\ldots,\lambda_{n}$ of the matrix $A$ are called the \textbf{multipliers} of fixed point.
We denote by $n_{0}$ the number of multipliers on the unit circle, by $n_{-}$ the number of multipliers inside the unit circle
and by $n_{+}$ the number of multipliers which lies outside the unit circle.  The fixed point is
\textbf{hyperbolic} if $n_{0}=0$, that is, there is no multipliers on the unit circle. Hyperbolic point is called a hyperbolic saddle if $n_{-} n_{+}\neq 0$. The fixed point is \textbf{nonhyperbolic} if $n_{0}\neq 0$.
We will also need the following definition (see \cite{laho1}).
\begin{defin}
Let $F:(x_0-r,x_0+r)\rightarrow \mathbb{R}$, $r>0$, be a map of class $C^{k}$, and $x_0$ is a fixed point of $F$ such that $F'(x_0)\neq 0$. 
If there is a $k\geq 3$ such that 
$F''(x_0)=\ldots=F^{(k-1)}(x_0)=0$ and $F^{(k)}(x_0)\neq 0$, then we say that the map $F$ is a \textbf{$k$-nondegenerate map in $x_0$}.
Specially, if $F''(x_0)\neq 0$, then we say that $F$ is \textbf{$2$-nondegenerate map in $x_0$}.
The number $k$ is called the \textbf{order of nondegeneracy} of map $F$ in $x_0$.
\end{defin}

In this paper the main object of our study is a box dimension of the orbit around the nonhyperbolic fixed point of discrete planar dynamical system. 
In Section 2, we prove the result for the two-dimensional systems with only one multiplier on the unit circle by using the center manifold theory. In Section 3, we analyze the Neimark-Sacker bifurcation. We will see that there is a difference in box dimension in rational and irrational case. The rational case will be proven by direct calculation, while the irrational case is showed using the analogous result for the Hopf bifurcation from the article \cite{zuzu}.

\section{Fractal analysis of bifurcations in $\mathbb{R}^2$ }

We will begin the study of bifurcation of discrete dynamical systems in $\eR^2$ with the bifurcations of two dimensional maps with only one multiplier on the unit circle.
Without loss of generality, the system can be written in the form 
\begin{eqnarray} \label{s0}
x_{n+1}&=&\lambda_1 x_{n} + f(x_{n},y_{n})\\
y_{n+1}&=&\lambda_2 y_{n} + g(x_{n},y_{n}) \nonumber
\end{eqnarray}
where $\left|\lambda_1\right|=1$ i $\left|\lambda_2\right|<1$, and $f$, $g$ are of class $C^{r}$ on some neighbourhood around the origin such that $f(0,0)=0$, $Df(0,0)=0$, $g(0,0)=0$, $Dg(0,0)=0$. It is clear that $(0,0)$ is a nonhyperbolic fixed point of the system. All other cases for multipliers can be proven analogously.
Then by Center Manifold Theorem we know that there exists $C^{r}$-center manifold $y=h(x)$, and the restriction of the system is
\begin{eqnarray} \label{s1}
x_{n+1}&=&\lambda_1 x_{n} + f(x_{n},h(x_{n}))\\
y_{n+1}&=&h(x_{n+1}). \nonumber
\end{eqnarray}

Now we would like to determine the connection between the box dimension of the restriction and its projection on the $x$-axis which is
\begin{eqnarray} \label{s2}
x_{n+1}&=&\lambda_1 x_{n} + f(x_{n},h(x_{n}))\\
y_{n+1}&=&0. \nonumber
\end{eqnarray}
We know that the projection on the $x$-axis is a Lipschitz map, but in order to get the equality of the dimensions we need the bilipschitz mapping. 
At this point, we will use the next lemma from the article \cite{zuzu3}.

\begin{lemma} (see \cite{zuzu3})
Let $g:\mathbb{R}^{N-1}\rightarrow \mathbb{R}$ be a Lipchitz map, $N\geq 2$, and we define $F:\mathbb{R}^{N}\rightarrow\mathbb{R}^{N}$ with $F(x,z):=(x,z+g(x))$, where $x\in\mathbb{R}^{N-1}$ and $z\in\mathbb{R}$. Then $F$ is a bilipschitz map and holds the measure, that is, for every measurable set $E\subset\mathbb{R}^{N}$ of limited measure it holds $\left|F(E)\right|=\left|E\right|$. Furthermore, for every limited set $A\subset\mathbb{R}^{N}$ we have
$$\overline{\dim}_{B}F(A)=\overline{\dim}_{B}A\,\,\, \rm{i}\,\,\, \underline{\dim}_{B}F(A)=\underline{\dim}_{B}A.$$
Set $A$ is a nondegenerate if and only if $F(A)$ is a nondegenerate.
\end{lemma}

We apply this lemma on the plane ($N=2$) and can see that if $g:\mathbb{R}\rightarrow\mathbb{R}$ is a Lipschitz map, then the map
$F:\mathbb{R}^2\rightarrow \mathbb{R}^2$ defined by
\begin{equation} \label{s3}
F(x,z)=(x,z+g(x))
\end{equation}
for $x$,$z\in\mathbb{R}$ is a bilipscitz map.\\

Now by using the above lemma we can prove the theorem about the box dimension of a discrete system on the center manifold.

\begin{theorem}
Let the restriction of the system ($\ref{s0}$) with $\lambda_1=\pm 1$ on the center manifold $y=h(x)$ is given by
$$x\mapsto G(x)=\lambda_1 x + f(x,h(x)).$$
Let $S(x_1,y_1)=(x_{n},y_{n})$ be an orbit of the systems on the center manifold in the form 
\begin{eqnarray} \label{centmnog}
x_{n+1}&=&\lambda_1 x_{n} + f(x_{n},h(x_{n}))\\
y_{n+1}&=&h(x_{n+1}) \nonumber
\end{eqnarray}
with initial point $(x_1,y_1)$ near $(0,0)$. If the map $G$ is a $k$-nondegenerate in $x=0$, 
then there exists $r>0$ such that for $|x_1|<r$ we have $$\dim_{B}S(x_1,y_1)=1-\frac{1}{k}.$$
\end{theorem}

\textbf{Proof.} First we define the set $A=\{(x,y):x\in A(x_1),y=0\}$, where $A(x_1)=(x_{n})_{n\in\mathbb{N}}$ is a one-dimensional discrete dynamical system generated by $x_{n+1}=G(x_{n})=\lambda_1 x_{n} + f(x_{n},h(x_{n}))$ and $x_1\in(0,r)$. Then we act with the map ($\ref{s3}$) on the set $A$, i.e. with  $F(x,z)=(x,z+h(x))$ and get $$F(A)=F(x,0)=(x,0+h(x))=(x,h(x)).$$ 
So the image of the set $A$ under map $F$ is associated system $(\ref{centmnog})$ on the center manifold $y=h(x)$. In other words, the map $F$ 
associate the projection of a orbit on $x$-axis with appropriate orbit on the center manifold. 
Since the map $h$ is of class $C^{r}$ on some neighbourhood small enough $\left|x\right|<\delta$ and $h'(0)=0$, then $h'$ is limited on that neighbourhood, and we have
$$\left|h(x_2)-h(x_1)\right|=\left|h'(x^{*})\right|\left|x_2-x_1\right|\leq M\left|x_2-x_1\right|$$
for some $x^{*}\in(x_1,x_2)$. Therefore $h$ is a Lipschitz map for $\left|x\right|<\delta$.
Now it follows from Lemma 1 that $F$ is a bilipschitz map, and 
$$\dim_{B} F(A)=\dim_{B}A.$$
Since $G$ is a $k$-nondegenerate map, from \cite{laho1} (Theorem 2.2), it follows that there exists $r>0$ such that for the sequence $A(x_1)=(x_{n})_{n\in\mathbb{N}}$ defined by $x_{n+1}=G(x_{n})$, $\left|x_1\right|<r$ we have
$$\dim_{B} A(x_1)=1-\frac{1}{k}.$$
We see now that $\dim_{B}A=\dim_{B}A(x_1)=1-\frac{1}{k}$.
Notice that $S(x_1,y_1)=F(A)$, so it follows that for $|x_1|<r$ is
$$\dim_{B}S(x_1,y_1)=\dim_{B} F(A)=\dim_{B}A=1-\frac{1}{k}.\,\,\,\,\blacksquare$$ 

\textbf{Remark 1.} The previous result can be easily generalized to the dynamical systems in $\mathbb{R}^{n}$ with only one multiplier on the unit circle.\\
 
Next we will apply the previous result to the systems with $\lambda_1=1$ or $\lambda_1=-1$. By short calculation we get the following forms of the restrictions and center manifolds.

The center manifold of the system ($\ref{s0}$) with $\lambda_1=\pm 1$ is given by
$$v=V(u)=\frac{1}{2}\omega u^2 + \mathcal{O}(u^3),$$
and the restriction of the system ($\ref{s0}$) on that center manifold is
\begin{equation} \label{s4}
x\mapsto \pm x + \frac{1}{2}\sigma x^2 + \frac{1}{6}(\delta + \frac{3ab}{1-c})x^3 + \mathcal{O}(x^4),
\end{equation}
where $\delta,\sigma,a,b\in\mathbb{R}$ are specified by 
$$\sigma=\frac{\partial^2 f(0,0)}{\partial x^2},\,\,\,\, \delta=\frac{\partial^3 f(0,0)}{\partial x^3},\,\,\,\,a=\frac{\partial^2 g(0,0)}{\partial x^2},\,\,\,b=\frac{\partial^2 f(0,0)}{\partial x \partial y}.$$

Notice that in the cases of saddle-node, pitchfork, transcritical and period doubling bifurcations the box dimension can be established from ($\ref{s4}$). For details see \cite{laho1}, \cite{kuz}, \cite{pe}.
So, in the case of saddle-node and transcritical bifurcation $\sigma\neq 0$, so the box dimension of a orbit of the system on the center manifold for $x_1$ small enough is
$$\dim_{B}S(x_1,y_1)=\frac{1}{2}.$$ 
For pitchfork and period doubling bifurcation is $a=b=0$ and $\delta\neq0$ so we have
$$\dim_{B}S(x_1,y_1)=\frac{2}{3}.$$
Now we can look at some examples.\\

\textbf{Example 1.} We consider the family of discrete planar systems of a form
\begin{eqnarray}
x &\mapsto& \mu + x + x^2\\
y &\mapsto& y/2 \nonumber
\end{eqnarray}
with parameter $\mu\in\mathbb{R}$.
Notice that $a=0$ so $\omega_1=0$ and the centre manifold is $v=0$. From $\sigma=2$ and $\delta=0$, it follows that the restriction of the system for  $\mu=0$ on the center manifold $v=0$ is of a form
$$u\mapsto u + u^2.$$ As we can see at Figure 1, the restriction of the system on the center manifold exhibits saddle-node bifurcation, and for $\mu=0$ the fixed point $(0,0)$ is nonhyperbolic with the box dimension $\dim_{B}S=\frac{1}{2}$. After the bifurcation we have two fixed points on the $x$-axis. Also notice that the box dimension of orbits of one-dimensional system $y\mapsto y/2$ on $y$-axis near $y_0=0$ is $0$, since $y_0=0$ is a hyperbolic fixed point for that map. So we can say that we have \textit{hyperbolic} and \textit{nonhyperbolic} direction.
\begin{center}
\epsfig{file=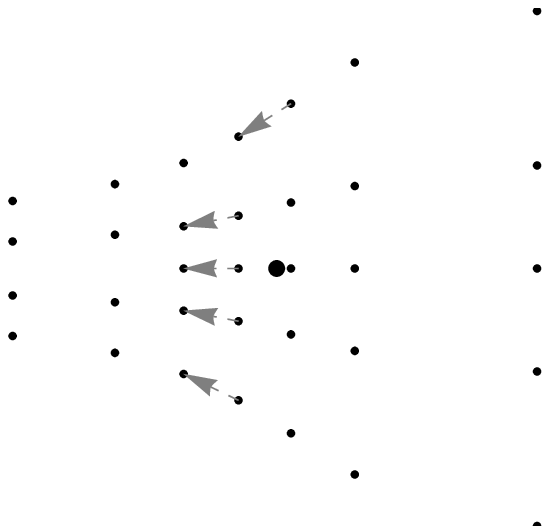, height=4cm}\epsfig{file=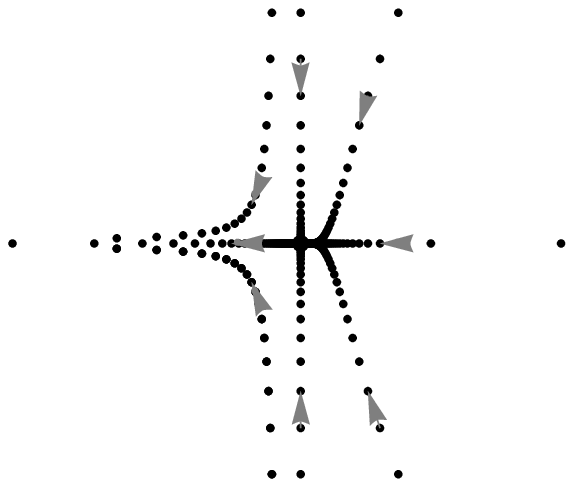, height=5cm}\\
\textit{\textbf{Figure 1a} $\mu<0$}\hskip 3cm\textit{\textbf{Figure 1b} $\mu=0$} \\
\epsfig{file=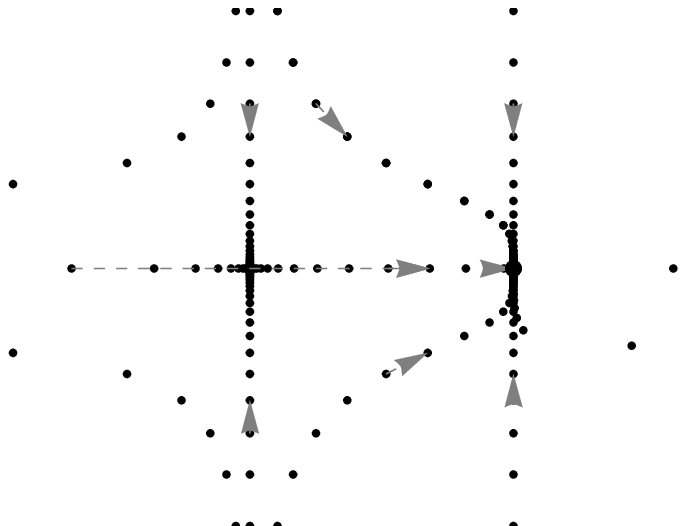, height=4cm}\\
\textit{\textbf{Figure 1c} $\mu>0$}
\end{center}
\vskip 1cm
\textbf{Example 2.} We consider the family
\begin{eqnarray}
x &\mapsto& (\mu-1) x + x^3\\
y &\mapsto& y/2 \nonumber
\end{eqnarray}
with the parameter $\mu\in\mathbb{R}$.
Notice that $a=0$ so $\omega_2=0$ and the center manifold is $v=0$. Also, notice that $\sigma=0$ so the restriction of the system for $\mu=0$ on the center manifold $v=0$ is of a form
$$u\mapsto -u + u^3.$$
As we can see at Figure 2, the restriction of the system on the center manifold exhibits the period doubling bifurcation, so for $\mu=0$ the fixed point $(0,0)$ is nonhyperbolic with box dimension $\dim_{B}S=\frac{2}{3}$. 
\begin{center}
\epsfig{file=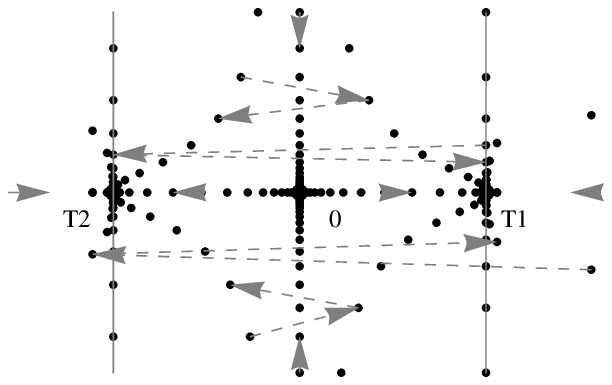, height=5cm}\\
\textit{\textbf{Figure 2a} $\mu<0$}\\
\epsfig{file=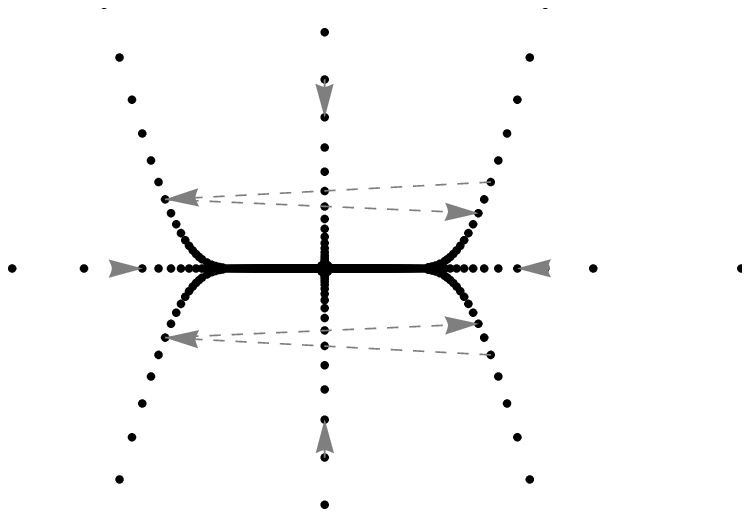, height=4cm}\epsfig{file=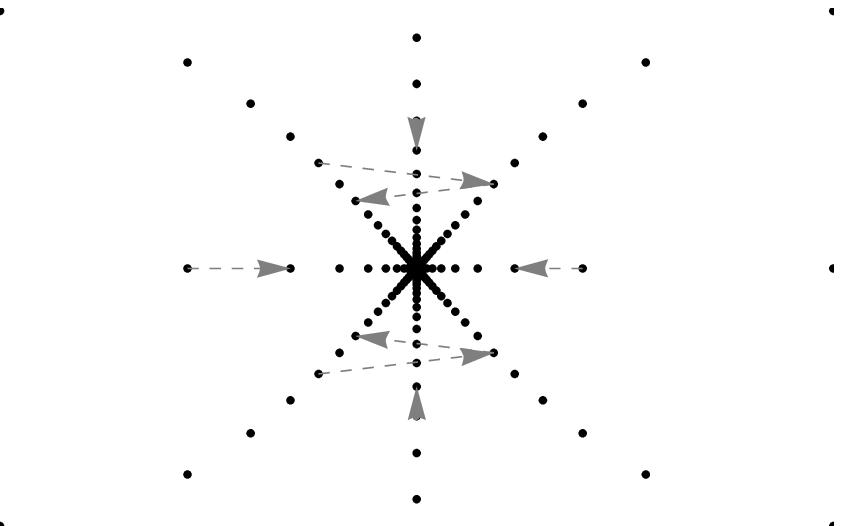, height=4cm}\\
\textit{\textbf{Figure 2b} $\mu=0$}\hskip 3cm\textit{\textbf{Figure 2c} $\mu>0$}
\end{center}

Considering the bifurcations with only one multiplier on the unit circle, we see that in this cases only more fixed or periodic points can bifurcate from one nonhyperbolic fixed point. In order to get the objects of bigger dimension from the bifurcation (e.g. invariant circle), we must have at least two multipliers on the unit circle.

The above-mentioned procedure for getting the box dimension of the orbit on the center manifold can be analogously applied to the orbits on the stable and unstable manifolds (see \cite{laho2}). This result can also be applied to appropriate bifurcations of continuous dynamical systems (see \cite{laho2}).

\section{Fractal Analysis of Neimark-Sacker bifurcation}

In this section we study the bifurcation at the nonhyperbolic fixed point with two complex conjugate eigenvalues on the unit circle when invariant curve bifurcates called Neimark-Sacker bifurcation. For details, see \cite{kuz}, \cite{ne},\cite{sa}.

From \cite{kuz}, we know that the polar normal form of one-parameter family of maps for Neimark-Sacker bifurcation is
\begin{eqnarray} \label{parfamns}
r&\mapsto& r + d\mu r + ar^3 + \mathcal{O}(\mu^2 r,\mu r^3,r^4)=F(r)\nonumber\\
\varphi &\mapsto& \varphi + \Theta_0 + \Theta_1 \mu + br^2 + \mathcal{O}(\mu^2,\mu r^2,r^3)=G(\varphi,r)
\end{eqnarray}
where $\mu \in \mathbb{R}$, and coefficients $d,a,\Theta_0,\Theta_1,b\in \mathbb{R}$ satisfy the bifurcation conditions: $a\neq 0$, $d\neq 0$ and $e^{ik\Theta_0}\neq 1$ for $k=1,2,3,4$.
Since the map $F$ doesn't depend on the angle $\varphi$, studying the existence of the invariant curve is simple. On the other hand, the map $G$ depends on $r$, which complicates the study of the orbit structure on the invariant curve, and consequently, the calculation of box dimension. All the previously studied bifurcations (see \cite{neveda}, \cite{zuzu}, \cite{laho1}) showed that the box dimension of a orbit around the nonhyperbolic fixed point is connected with the box dimension of the invariant set which emerge at the bifurcation point. For instance, in the bifurcation of one-dimensional discrete dynamical systems when only fixed point can bifurcate, the box dimension is between 0 and 1. Furthermore, the box dimension in Hopf bifurcation when we have limit cycles is greater  then 1.

Moreover, in all other bifurcations of discrete and continuous systems the invariant sets which emerge consist of only one orbit. For example, the limit cycle which bifurcates from the Hopf bifurcation has box dimension 1, and consist of only one trajectory.  In the case of Neimark-Sacker bifurcation, the originate invariant curve consists of many different orbits, and their structures depend on the map $G(\varphi,r_0)$ from ($\ref{parfamns}$), with $r=r_0$ the radius of invariant curve. 

We know that the orbit structure on the invariant curve is depended only on the map $G(\varphi,r)$. Therefore, there is a case 
$G(\varphi,r_0)=\varphi+\Theta_0$, $\Theta_0=\frac{2\pi p}{q}$ $\rightarrow$. Then $(p,q)$-cycles emerge by the bifurcation and the box dimension of every orbit on the invariant curve is 0 (property of finite stability).
In the case $G(\varphi,r_0)=\varphi+\Theta_0$, $\Theta_0=2\pi\beta$, $\beta$ irrational, on the invariant curve every orbit is dense and its box dimension is $1$. It is known that in general case with $G(\varphi,r_0)=\varphi+\Theta_0+\Theta_1\mu + br_0^2 + \mathcal{O}(\mu^2,\mu r_0^2,r_0^3)$ $\rightarrow$, the orbit structure changes by changing the parameter between rational and irrational rotation numbers (Arnold tongues, see \cite{kuz},\cite{ap}).

In the article \cite{laho2}, it was showed that in the neighbourhood of the (un)stable hyperbolic fixed point in $\eR^{n}$, the box dimension of every orbit is 0. Now we want to establish what happens with the box dimension of an orbit around the nonhyperbolic fixed point.
So there is no need to study whole family of one-parameter maps but only the system for the bifurcation value $\mu=0$ because then the point $x_0=(0,0)$ is a nonhyperbolic fixed point, and the box dimension is positive. Our goal is to determine the value of that box dimension. 

Hence, we look at the polar normal form for Neimark-Sacker bifrucation at the bifurcation value $\mu=0$ 
\begin{eqnarray} \label{polnf}
r&\mapsto& r + ar^3 + \mathcal{O}(r^4)=f(r) \nonumber \\ 
\varphi &\mapsto& \varphi + \Theta_0 + br^2 + \mathcal{O}(r^3)=g(\varphi,r)
\end{eqnarray}
where $a\neq 0$ (nondegeneracy condition) and $e^{i\Theta}\neq\sqrt[n]{1}$, for $n=1,2,3,4$ (nonresonant condition).
We consider the orbits of the system ($\ref{polnf}$) around the nonhyperbolic fixed point $x_0=(0,0)$. 
Let $a<0$, and define
\begin{equation} \label{disspirala}
\Gamma=\{(r_{k},\varphi_{k}): r_{k}=f(r_{k-1}),\varphi_{k}=g(r_{k-1},\varphi_{k-1}),k\in\mathbb{N}, (r_0,\varphi_0) \,\,\rm{given} \}.
\end{equation}
When we draw only one orbit $\Gamma$, we observe that the points are spirally going to the origin, and that is why the orbit $\Gamma$ is called a \textbf{discrete spiral}. To be precise, for $a<0$ the discrete spiral is spirally going to the origin, but for $a>0$ the orbit is spirally going away from the origin. Furthermore, if the rotation angle is positive, then the spiral has a positive direction, and otherwise the direction is negative. See Figure 3.

\begin{center}
\includegraphics[width=5.5cm]{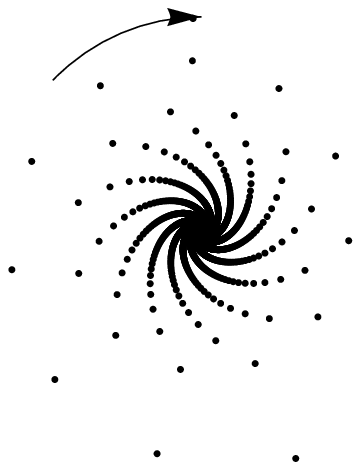}  \includegraphics[width=5.5cm]{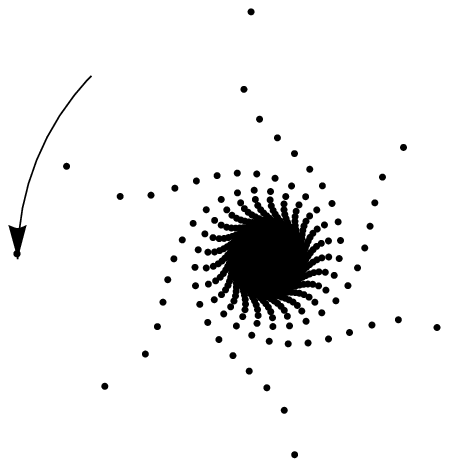}\\ 
\textit{\textbf{Figure 3a} $f(r)=r-r^3$,} \hskip 3cm \textit{\textbf{Figure 3b} $f(r)=r-r^3$},\\ 
\textit{ $g(r,\varphi)=\varphi-\frac{\pi}{6}-r^2$}\hskip 5cm \textit{$g(r,\varphi)=\varphi+1+r^2$} \\
\end{center}

\subsection{Bounds for box dimension}

We know that for every bounded set in the plane we have
$$0\leq\underline{\dim}_{B}\Gamma\leq\overline{\dim}_{B}\Gamma\leq 2.$$ 

Now we recall that the generalized polar normal form for the Neimark-Sacker bifurcation proven in the article \cite{bave1} (Theorem 5) is
\begin{eqnarray} \label{genpolnf}
r&\mapsto& r -\sum_{m=1}^{p} a_{2m+1}r^{2m+1} + \mathcal{O}(r^{2p+2})=f(r) \nonumber \\ 
\varphi &\mapsto& \varphi + \Theta_0 + \sum_{m=1}^{p}b_{2m}r^{2m} + \mathcal{O}(r^{2p+1})=g(\varphi,r).
\end{eqnarray}
and the corresponding nondegeneracy condition:
\begin{equation}
a_3=\ldots=a_{2k-1}=0,\,\,\,a_{2k+1}\neq 0.
\end{equation}
It is clear that $\alpha=2k+1$ is the order of nondegeneracy of the map $f$.

\begin{theorem} \textbf{Bounds for box dimension}\\
Let $\Gamma$ be a discrete spiral defined by ($\ref{disspirala}$) where $f$ and $g$ are defined by the system $(\ref{genpolnf})$. If $f$ is an $\alpha$-nondegenerate in $x_0=0$, then there exist $r_0>0$ small enough such that
$$1-\frac{1}{\alpha}\leq\underline{\dim}_{B}\Gamma(r_0,\varphi_0)\leq\overline{\dim}_{B}\Gamma(r_0,\varphi_0)\leq 2(1-\frac{1}{\alpha}).$$
\end{theorem}

\textbf{Proof.}\\
If $f$ is an $\alpha$-nondegenerate map in a point $x_0$, it means that
$$f''(x_0)=\ldots=f^{(\alpha-1)}(x_0)=0,\,\,f^{(\alpha)}(x_0)\neq 0.$$
In the case when $f=f(r)$ in the system $(\ref{genpolnf})$, then
$$a_{3}=\ldots=a_{2k-1}=0,\,\,a_{2k+1}\neq 0$$
for some $k\in{1,2,\ldots,p}$, so then $\alpha=2k+1$.\\
Let $A$ be a sequence of concentric circles of radius $r_{k+1}=f(r_{k})$, $k\in \mathbb{N}$ with the centre in the origin. Notice that $\Gamma\subset A$, so in the case of nonhyperbolic fixed point
$$\overline{\dim}_{B}\Gamma \leq \dim_{B} A=2\dim_{B} \{r_{k}\}_{k\in\mathbb{N}}=2(1-\frac{1}{\alpha}).$$
Of course, we have used the fact that the box dimension of the orbit of one-dimensional sequence $\{r_{k}\}$ generated by an $\alpha$-nondegenerate map is $1-\frac{1}{\alpha}$ (see \cite{laho1}). In this way we got the upper bound for the box dimension of $\Gamma$ which is, in fact, the best upper bound.\\
Now we are interested in the lower bound.
Let $P_{x}:\mathbb{R}^2\rightarrow \mathbb{R}^2$ be a radial projection on the $x$-axis defined by $P_{x}(r,\varphi)=(r,0)$. 
This map is a Lipschitz map, so 
$$\underline{\dim}_{B}\Gamma\geq \dim_{B} P_{x}(\Gamma)=\dim_{B} \{r_{k}\}_{k\in\mathbb{N}}.$$

It means that we get the following estimation which holds around the nonhyperbolic fixed point
\begin{eqnarray}
\dim_{B} \{r_{k}\}_{k\in\mathbb{N}}\leq \underline{\dim}_{B} \Gamma \leq\overline{\dim}_{B}\Gamma \leq 2\dim_{B} \{r_{k}\}_{k\in\mathbb{N}},\nonumber
\end{eqnarray}
that is,
\begin{equation}
1-\frac{1}{\alpha}\leq\underline{\dim}_{B}\Gamma\leq\overline{\dim}_{B}\Gamma\leq 2(1-\frac{1}{\alpha}).\,\,\,\blacksquare \nonumber
\end{equation}
So, in the case of classic Neimark-Sacker bifurcation for $\alpha=3$ we have
\begin{eqnarray}
\frac{2}{3}\leq \underline{\dim}_{B} \Gamma\leq\overline{\dim}_{B}\Gamma \leq \frac{4}{3}.
\end{eqnarray}

In the next section we will shown that for the fixed and rational rotation angle the lower bound is achieved, while in the irrational case, the upper bound is achieved.

\subsection{Fixed and rational angle displacement map}

The polar normal form of Neimark-Sacker bifurcation at the bifurcation value $\mu=0$ with fixed rotation angle ($g$ doesn't depend on $r$) 
\begin{eqnarray} \label{genpolnf2}
r &\mapsto& r + ar^{\alpha} + \mathcal{O}(r^{\alpha+1})=f(r) \nonumber \\ 
\varphi &\mapsto& \varphi + \Theta_0 =g(\varphi)
\end{eqnarray}
where $\alpha$ is odd. If the rotation angle is of a form
\begin{equation} \label{racrot}
\Theta_0=\frac{2\pi p}{q},
\end{equation}
for some $p,q \in \mathbf{Z}$, $N(p,q)=1$, we say that the rotation angle is rational. That is the special case of Neimark-Sacker bifurcation, when by bifurcation emerges the invariant circle with the $(p,q)$-periodic cycles on it, that is, every orbit on the invariant circle consist of finite number of points, so its box dimension is 0.
Due to the fact that there exists $q$ such that
$$q\Theta_0=2\pi p,$$ the orbit $\Gamma$ of the system ($\ref{genpolnf2}$) will be the union of $q$ rays from the origin with the one-dimensional sequences of points which tend to origin. On Figure 4 we can see the phase portraits of the system 
\begin{eqnarray}
f(r)&=&(1+\mu)r-r^3\nonumber\\
g(r,\varphi)&=&\varphi+\frac{\pi}{6}
\end{eqnarray}
depending on $\mu$. Notice that $q=12$, so on every picture the points are on the rays, but around the nonhyperbolic fixed point they are evident (Figure 4b). On Figure 4c  we see that the orbit on the invariant circle are $(12)$-cycles.

\begin{center}
\includegraphics[width=4.2cm]{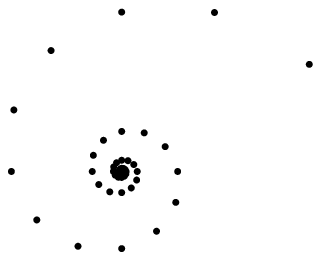}\includegraphics[width=4.2cm]{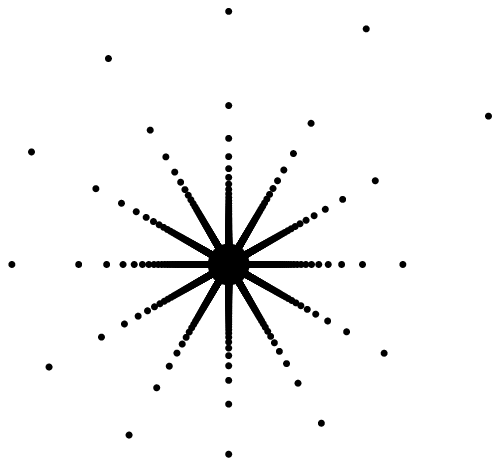}\includegraphics[width=4.2cm]{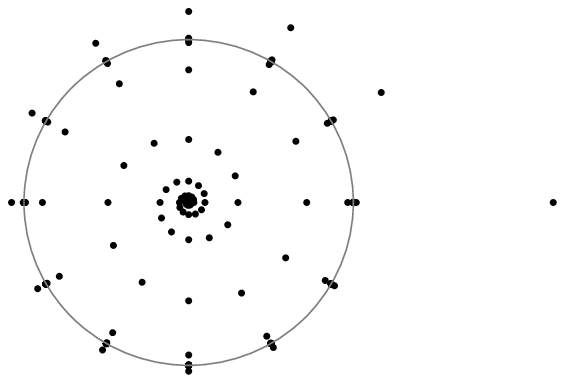}\\ 
\textit{\textbf{Figure 4a} $\mu<0$} \hskip 1.2cm \textit{\textbf{Figure 4b} $\mu=0$} \hskip 2cm \textit{\textbf{Figure 4c} $\mu>0$}
\end{center}

For the box dimension result we will need following lemma.
\begin{lemma} \label{kiteracijaf}
Let be a system
\begin{eqnarray}
r &\mapsto& r + ar^{\alpha} + \mathcal{O}(r^{\alpha+1})=f(r) \nonumber \\ 
\varphi &\mapsto& \varphi + \Theta_0 +br^{\alpha-1} + \mathcal{O}(r^{\alpha})=g(\varphi,r)
\end{eqnarray}
with $\alpha\geq 3$.
Then the $k.$ iteration of that system is of a form
\begin{eqnarray}
f^k(r)&=&r+kar^{\alpha}+\mathcal{O}(r^{\alpha+1})\nonumber\\
g^{k}(\varphi,r)&=&\varphi+k\Theta_0 + bkr^{\alpha-1}+\mathcal{O}(r^{\alpha}).
\end{eqnarray}
\end{lemma}

\textbf{Remark 2.} This technical lemma can be easily proven by using the mathematical induction.

\begin{theorem} (\textbf{fixed and rational rotation angle })\\
Let $\Gamma(r_0,\varphi_0)$ be a discrete spiral defined by ($\ref{disspirala}$) where $f$ and $g$ are given by a system $(\ref{genpolnf2})$, 
so let the fixed rotation angle be in a form ($\ref{racrot}$). Then for $r_0>0$ small enough we have
$$\dim_{B}\Gamma(r_0,\varphi_0)=1-\frac{1}{\alpha}.$$
Moreover, $\Gamma(r_0,\varphi_0)$ is Minkowski nondegenerate.
\end{theorem}

\textbf{Proof.}\\
Because of the fixed and rational rotation angle, for the discrete spiral with the initial point $(r_0,\varphi_0)$, it holds
\begin{eqnarray}
\Gamma(r_0,\varphi_0)=\cup_{i=0}^{q-1} \Gamma_{i}(r_0,\varphi_{i})
\end{eqnarray}
where
\begin{eqnarray} \Gamma_{i}(r_0,\varphi_{i})&=&\{f^{i}(r_0),f^{q+i}(r_0),f^{2q+i}(r_0),\ldots\}\nonumber\\
&=&\{f^{i}(r_0),f^{q}(f^{i}(r_0)),f^{2q}(f^{i}(r_0)),\ldots\}
\end{eqnarray}
is a one-dimensional sequence on the ray $\varphi=\varphi_{i}=g^{i}(\varphi_0)$. 
We see that the discrete spiral $\Gamma_{i}(r_0,\varphi_{i})$ is in fact generated by the map
 $f^{q}$ with the initial point $f^{i}(r_0)$, $i=0,\ldots,q-1$. 
From Lemma 2, we have $f^{q}(r)=r+qar^{\alpha}+\mathcal{O}(r^{\alpha+1})$. 
 
Now from \cite{laho1} (Theorem 2.2), it follows that for $r_0>0$ small enough
$$\dim_{B}\Gamma_{i}(r_0,\varphi_{i})=1-\frac{1}{\alpha}.$$
Now by using the finite stability property (see Section 1) we get
\begin{equation}
\dim_{B}\Gamma(r_0,\varphi_0)=\dim_{B}\Gamma_{1}(r_0,\varphi_0)=1-\frac{1}{\alpha}
\end{equation}
since all the box dimension on the rays are equal for $r_0$ small enough. \\

Regarding the Minkowski content, it can be easily seen that
$$\left|\Gamma_{1,\varepsilon}\right| \leq\left|\Gamma_{\varepsilon}\right| \leq \sum_{i=1}^{q-1} \left|\Gamma_{i,\varepsilon}\right|,$$
so
$$\mathcal{M}_{*}^{d}(\Gamma_1)\leq \mathcal{M}_{*}^{d}(\Gamma)\leq\mathcal{M}^{*d}(\Gamma)\leq \sum_{i=1}^{q-1} \mathcal{M}^{*d} (\Gamma_{i})$$
where $d=\dim_{B}\Gamma=\dim_{B} \Gamma_{i}$. From \cite{laho1} (Theorem 2.2) we know that $\Gamma_{i}$ are Minkowski nondegenerate sets, so we get that
$$0<\mathcal{M}_{*}^{d}(\Gamma)\leq\mathcal{M}^{*d}(\Gamma)<\infty.\,\,\,\,\blacksquare$$

\subsection{General rational case}

In this section we will show the method for getting the box dimension in general rational case. In the previous section we got the box dimension result by using the finite stability property and some previously known results. In this section, we will show the method of direct calculation by using the definition of box dimension. First step is finding the good estimation of $\left|\Gamma_{\varepsilon}\right|$, where $\Gamma=\Gamma(r_0,\varphi_0)$ is a discrete spiral and $\Gamma_{\varepsilon}$ is its $\varepsilon$-neighbourhood. 

At the beginning, we will study the calculation method with the fixed rotation number.
Hence, we would like to estimate the area of the $\varepsilon$-neighbourhood of one orbit of the system ($\ref{genpolnf2}$)
with $a<0$, fixed and rational rotation angle  $\Theta_0$ and for fixed $\varepsilon>0$. Without loss of generality we may assume that $\alpha=3$. 
For odd $\alpha>3$ the calculus is analogous. In order to do that, we have to look in detail what are all the possible overlaps in the mentioned $\varepsilon$-neighbourhood. Notice that the following analysis is valid in the general situation as well. 

In the case of discrete orbit $\Gamma$, the $\varepsilon$-neighbourhood of $\Gamma$, we denote it by $\Gamma_{\varepsilon}$, is a union of circles around points $A_{k}(r_{k},\varphi_{k})$ with radius $\varepsilon$. The main problem is that after some $k$ the circles begin to overlap and then it is difficult to calculate its overall surface. 
We observe two neighbouring points $A_{k}(r_{k},\varphi_{k})$ and $A_{k+1}(r_{k+1},\varphi_{k+1})$, and one point in the next level $A_{k+q_0}(r_{k+q_0},\varphi_{k+q_0})$ where $q_0$ is the lowest positive number such that $q_0\Theta_0>2\pi$. See Figure 5. We will see that in the case of rational rotation number with $p=1$ is $q_0\Theta_0=2\pi$. For $p>1$ the number $q_0$ will be such that $q_0\Theta_0=2p\pi$.
Therefore, we have three significant distances whose behavior we must analyze in order to determine the way and the order of overlaps.

\begin{center}
\includegraphics[width=5cm]{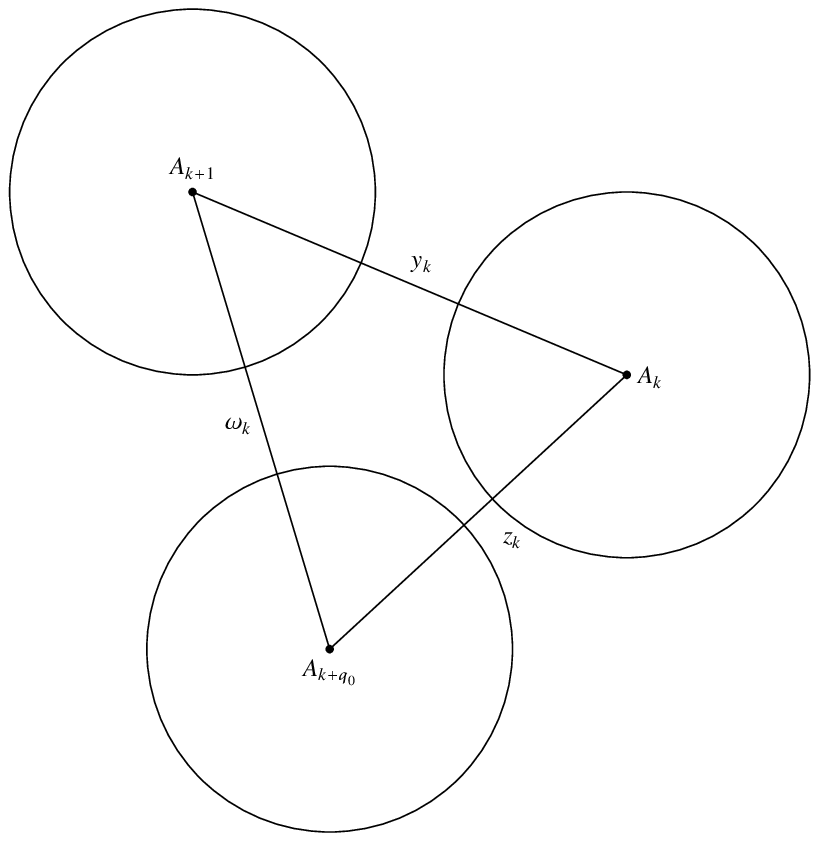}\\ 
\textit{\textbf{Figure 5}} 
\end{center}

We denote with
\begin{eqnarray} \label{ponnizy}
y_{k}&=&d(A_{k},A_{k+1})=\sqrt{r^2_{k}+r^2_{k+1}-2r_{k}r_{k+1}\cos(\Theta_0)}
\end{eqnarray}
\begin{eqnarray} \label{ponnizz}
z_{k}&=&d(A_{k},A_{k+q_0})=\sqrt{r^2_{k}+r^2_{k+q_0}-2r_{k}r_{k+q_0}\cos(q_0\Theta_0)}\\\label{ponnizw}
w_{k}&=&d(A_{k+1},A_{k+q_0})=\sqrt{r^2_{k+1}+r^2_{k+q_0}-2r_{k+1}r_{k+q_0}\cos((q_0-1)\Theta_0)}.
\end{eqnarray}

Hence, we have three possibilities 
\begin{eqnarray}
y_{k}&<&2\varepsilon \\
z_{k}&<&2\varepsilon\\
w_{k}&<&2\varepsilon
\end{eqnarray} 
First, for every overlapping we must show that the sequences $\{y_{k}\}$, $\{z_{k}\}$ and $\{\omega_{k}\}$ are decreasing to ensure that 
$m(\varepsilon)$ is the lowest natural number for which overlapping occurs. So we study the behaviour of the sequences
 $\{y_{k}\}$, $\{z_{k}\}$ and $\{w_{k}\}$ as $k\rightarrow \infty$.
 Since $r_{k+q_0}$ appears in the expressions for $z_{k}$ and $w_{k}$, we need the expression for $r_{k+q_0}=f^{q_0}(r_{k})$ showed in Lemma \ref{kiteracijaf}.

It is easy to prove the lemma about the behaviour of $\{y_{k}\}$, $\{z_{k}\}$ and $\{w_{k}\}$ for fixed and rational $\Theta_0$.
\begin{lemma}
Let ($\ref{genpolnf2}$) be a system with $a<0$. Then the sequences $\{y_{k}\}$, $\{z_{k}\}$ i $\{w_{k}\}$ defined by ($\ref{ponnizy}$), ($\ref{ponnizz}$) and ($\ref{ponnizw}$) are decreasing for $k$ big enough, and
\begin{eqnarray}
z_{k}&\simeq& q_0\left|a\right|k^{-\frac{\alpha}{\alpha-1}}\\
w_{k}\simeq y_{k}&\simeq& \sqrt{2(1-\cos(\Theta_0))}k^{-\frac{1}{\alpha-1}}.
\end{eqnarray}
Moreover, there exists $K_0$ such that $z_{k}<y_{k}$ for all $k>K_0$.
\end{lemma}

\textbf{Proof}.\\
First, we show that $\{y_{k}\}$ is decreasing ($y_{k+1}<y_{k}$):\\
By squaring $y_{k}$ we get 
$$r^2_{k+1}+r^2_{k+2}-2r_{k+1}r_{k+2}\cos(\Theta_0)<r^2_{k}+r^2_{k+1}-2r_{k}r_{k+1}\cos(\Theta_0)$$
and
$$r_{k}+r_{k+2}>2r_{k+1}\cos(\Theta_0).$$
In the proof of Theorem 2.2, \cite{laho1} it was shown that one-dimensional sequence of a form $r_{k+1}=r_{k}+ar_{k}^{\alpha}+\mathcal{O}(r_{k}^{\alpha+1})$, for $a<0$ has a decreasing sequence of differences, that is, $(r_{k}-r_{k+1})$ is decreasing on some small neighbourhood of the fixed point (for $k$ big enough). Since the box dimension is also calculated on some small neighbourhood of the fixed point, it means that for $k$ big enough it holds 
$$r_{k+1}-r_{k+2}<r_{k}-r_{k+1}$$
and follows
$$r_{k}+r_{k+2}>2r_{k+1}>2r_{k+1}\cos(\Theta_0).$$
So $\{y_{k}\}$ is decreasing for $k$ big enough.
Completely analogously can be proven that $\{z_{k}\}$ and $\{w_{k}\}$ are decreasing as well.\\

By including $r_{k+1}=r_{k}+ar_{k}^{\alpha} + \mathcal{O}(r_{k}^{\alpha+1})$ and $r_{k+q_0}=r_{k}+q_0ar_{k}^{\alpha}+\mathcal{O}(r_{k}^{\alpha+1})$ in the expressions ($\ref{ponnizy}$), $(\ref{ponnizz})$ and $(\ref{ponnizw})$ we get
\begin{eqnarray}
y_{k}^2&=&2(1-\cos(\Theta_0))r_{k}^2 + \mathcal{O}(r_{k}^{\alpha+1})\\
z_{k}^2&=&2(1-\cos(q_0\Theta_0))r_{k}^2 + \mathcal{O}(r_{k}^{\alpha+1})\\
w_{k}^2&=&2(1-\cos((q_0-1)\Theta_0))r_{k}^2 + \mathcal{O}(r_{k}^{\alpha+1})
\end{eqnarray}
 
Using the known behaviour of $\{r_{k}\}$ as $k\rightarrow \infty$ (Theorem 2.2, \cite{laho1}) 
$$r_{k}\simeq k^{-\frac{1}{\alpha-1}}\,\,\,\textrm{i.e.,}\,\,\,\exists A,B>0 \,\,\,\,Ak^{-\frac{1}{\alpha-1}}\leq r_{k} \leq Bk^{-\frac{1}{\alpha-1}} $$
we obtain
\begin{eqnarray}
y_{k} &\simeq& \sqrt{2(1-\cos(\Theta_0))}k^{-\frac{1}{\alpha-1}}\\ 
z_{k} &\simeq& \sqrt{2(1-\cos(q_0\Theta_0))}k^{-\frac{1}{\alpha-1}}\\
w_{k} &\simeq& \sqrt{2(1-\cos((q_0-1)\Theta_0))}k^{-\frac{1}{\alpha-1}}.
\end{eqnarray}
The condition $\cos(\Theta_0)\neq 1$ is satisfied because of the nonresonant condition for the Neimark-Sacker bifurcation.
This type of behaviours will be valid also for $\Theta_0=2\pi\beta$, $\beta$ irrational. 
Of course, if the rotation angle is rational, then $q_0\Theta_0=2\pi p$. Then it follows
\begin{eqnarray}
z_{k}&\simeq& q_0\left|a\right|k^{-\frac{\alpha}{\alpha-1}}\\
w_{k}\simeq y_{k}&\simeq& \sqrt{2(1-\cos(\Theta_0)}k^{-\frac{1}{\alpha-1}}.
\end{eqnarray}
Regarding the relation between $z_{k}$ and $y_{k}$, some minor problem can occur since both depend on $q_0$.
That is why we will closely examine their relationship.
Let $A_{z},B_{z}>0$ be constants such that
$$A_{z}q_0 \left|a\right|k^{-\frac{\alpha}{\alpha-1}}\leq z_{k} \leq B_{z}q_0\left|a\right|k^{-\frac{\alpha}{\alpha-1}},$$
while $A_{y},B_{y}>0$ are such that
$$A_{y}\sqrt{2(1-\cos(\Theta_0))}k^{-\frac{1}{\alpha-1}}\leq y_{k} \leq B_{y}\sqrt{2(1-\cos(\Theta_0))}k^{-\frac{1}{\alpha-1}}.$$
We want to have 
$$z_{k}\leq B_{z} q_0 \left|a\right|k^{-\frac{\alpha}{\alpha-1}}< A_{y}  \sqrt{2(1-\cos(\Theta_0)}k^{-\frac{1}{\alpha-1}}\leq y_{k},$$
that is,
$$k>\frac{B_{z}q_0 \left|a\right|}{A_{y}\sqrt{2(1-\cos\Theta_0)}}.$$
In order to obtain $\frac{1}{1-\cos(\Theta_0)}$ we use $\Theta_0=\frac{2\pi p}{q_0}$ and the Taylor expansion for cosine, and get the existence of  constants $C_1,C_2>0$ such that
$$\frac{q_0^2}{2\pi^2p^2C_1}\geq \frac{1}{1-\cos\Theta_0}\geq \frac{q_0^2}{2C_2\pi^2p^2}.$$
It means that $k$ need to satisfy
$$k>\frac{B_{z}\left|a\right|q_0^2}{2\pi p \sqrt{C_1}A_{y}}.$$
Now we have that $z_{k}<y_{k}$ for $k>K_0=\frac{B_{z}\left|a\right|q_0^2}{2\pi p\sqrt{C_1}A_{y}}$. $\blacksquare$\\

\textbf{Remark 3.} For $a>0$, we need to observe the inverse map for $f$, i.e. the system
\begin{eqnarray}
r &\mapsto& r - ar^{\alpha} + \mathcal{O}(r^{\alpha+1})\nonumber\\
\varphi &\mapsto& \varphi - \Theta_0
\end{eqnarray}
and the result from lemma also holds.\\

Now we study the case $\Theta_0=2\pi\beta$, $\beta\in\mathbb{Q}$, but the rotation angle isn't fixed but has also the terms of higher order.
So we observe the system
\begin{eqnarray} \label{opcrac}
r&\mapsto&r+ar^{\alpha}+\mathcal{O}(r^{\alpha+1})\nonumber\\
\varphi&\mapsto&\varphi + \Theta_0 + br^{\alpha-1} + \mathcal{O}(r^{\alpha})
\end{eqnarray}
where $\Theta_0=\frac{2\pi p}{q}$, $p,q\in\mathbb{Z}$ , $N(p,q)=1$, $a<0$ and $\alpha$ odd.

The overlapping sequences are given by 
{\small\begin{eqnarray} \label{ponorac1}
y_{k}&=&d(A_{k},A_{k+1})=\sqrt{r^2_{k}+r^2_{k+1}-2r_{k}r_{k+1}\cos(\Theta_0)}
\end{eqnarray}}
{\small\begin{eqnarray} \label{ponorac2}
z_{k}&=&d(A_{k},A_{k+q_0})=\sqrt{r^2_{k}+r^2_{k+q_0}-2r_{k}r_{k+q_0}\cos(\varphi_{k+q_0}-\varphi_{k})}
\end{eqnarray}}
{\small\begin{eqnarray} \label{ponorac3}
w_{k}&=&d(A_{k+1},A_{k+q_0})=\sqrt{r^2_{k+1}+r^2_{k+q_0}-2r_{k+1}r_{k+q_0}\cos(\varphi_{k+q_0}-\varphi_{k+1})}.
\end{eqnarray}}

\begin{lemma}
Let ($\ref{opcrac}$) be a system with $\alpha=3$. Then the sequences $\{y_{k}\}$, $\{z_{k}\}$ and $\{w_{k}\}$ defined by ($\ref{ponorac1}$), $(\ref{ponorac2})$ and $(\ref{ponorac3})$ are decreasing for $k$ big enough, and
\begin{eqnarray}
z_{k}&\simeq& q_0\sqrt{a^2+b^2}k^{-\frac{3}{2}}\\
w_{k}\simeq y_{k}&\simeq& \sqrt{2(1-\cos(\Theta_0))}k^{-\frac{1}{2}}.
\end{eqnarray}
Moreover, there exists $K_0>0$ such that $z_{k}<y_{k}$ for all $k>K_0$.
\end{lemma}

\textbf{Proof.}\\
First we will show that $\{z_{k}\}$ is decreasing for $k$ large enough. 
Hence, we have
$$z_{k}^2=r^2_{k}+r^2_{k+q_0}-2r_{k}r_{k+q_0}\cos(\varphi_{k+q_0}-\varphi_{k}),$$
and we want to show that $z_{k+1}<z_{k}$ for large $k$. 
It is easy to show that $\varphi_{k+q_0+1}-\varphi_{k+1}<\varphi_{k+q_0}-\varphi_{k}$ and it follows that
\begin{equation} \label{izrazcos}
-\cos(\varphi_{k+q_0+1}-\varphi_{k+1})<-\cos(\varphi_{k+q_0}-\varphi_{k}).
\end{equation}

Now recall that the sequence of differences for $\{r_{k}\}$ is also decreasing so 
$$r_{k}-r_{k+1}>r_{k+q_0}-r_{k+q_0+1}.$$
By squaring and reducing we get 
\begin{eqnarray}
r_{k}^2-r^{2}_{k+1}+r^2_{k+q_0}-r^2_{k+q_0+1}&>&2r_{k}r_{k+q_0}-2r_{k+1}r_{k+q_0+1}>\nonumber\\
&>&(2r_{k}r_{k+q_0}-2r_{k+1}r_{k+q_0+1})\cos(\varphi_{k+q_0}-\varphi_{k}),\nonumber
\end{eqnarray}
and by using ($\ref{izrazcos}$) we get
$$r_{k}^2+r_{k+q_0}^2-2r_{k}r_{k+q_0}\cos(\varphi_{k+q_0}-\varphi_{k})>r^2_{k+1}+r^2_{k+q_0+1}-2r_{k+1}r_{k+q_0+1}\cos(\varphi_{k+q_0+1}-\varphi_{k+1}).$$
Analogously for $\{w_{k}\}$ and $\{y_{k}\}$. 

Now we observe the behaviour of $z_{k}$. In the expression 
$$z_{k}^2=r^2_{k}+r^2_{k+q_0}-2r_{k}r_{k+q_0}\cos(\varphi_{k+q_0}-\varphi_{k})$$
we put
$$f^{q_0}(r)=r + q_0ar^3 + cr^4 + dr^5 +er^6 + \mathcal{O}(r^7)$$ and
$$\varphi_{k+q_0}=\varphi_{k}+q_0\Theta_0 + q_0br_{k}^2 + \mathcal{O}(r_{k}^4)$$ which we proved in Lemma $\ref{kiteracijaf}$
and
$$\cos(\varphi_{k+q_0}-\varphi_{k})=1-\frac{1}{2}b^2q_0^2r_{k}^4+\mathcal{O}(r_{k}^6)$$
and get
\begin{equation}
z_{k}^2=(a^2+b^2)q_0^2r_{k}^6 + \mathcal{O}(r_{k}^7).
\end{equation}

It is easy to show that the behaviour of $y_{k}$ is the same as for fixed and rational rotation angle (Lemma 3). 
The sequence $y_{k}$ is
\begin{eqnarray}
y^2_{k}&=& r_{k}^2 + r_{k+1}^2 - 2r_{k}r_{k+1}\cos(\varphi_{k+1}-\varphi_{k})=\nonumber\\
&=& r_{k}^2 + r_{k+1}^2 -2r_{k}r_{k+1}\cos(\Theta_0+br_{k}^2 + \mathcal{O}(r_{k}^4))
\end{eqnarray}
and the behaviour of cosine of the rotation angle is 
\begin{eqnarray}
\cos(\Theta_0+br_{k}^2 + \mathcal{O}(r_{k}^4))&=& \cos(\Theta_0)-(br_{k}^2 + \mathcal{O}(r_{k}^4))\sin(\Theta_0)-\nonumber\\
&-&(br_{k}^2 + \mathcal{O}(r_{k}^4))^2\frac{\cos(\Theta_0)}{2} + \mathcal{O}(r_{k}^6)=\nonumber\\
&=&\cos(\Theta_0)-br_{k}^2\sin(\Theta_0)+\mathcal{O}(r_{k}^4).
\end{eqnarray}
It follows
\begin{eqnarray}
y_{k}^2&=&r_{k}^2+(r_{k}+ar_{k}^3+\mathcal{O}(r_{k}^4))^2-\nonumber\\
&-&2r_{k}(r_{k}+ar_{k}^3+\mathcal{O}(r_{k}^4))(\cos(\Theta_0)-br_{k}^2\sin(\Theta_0)+\mathcal{O}(r_{k}^4))=\nonumber\\
&=& 2r_{k}^2 + 2ar_{k}^4+\mathcal{O}(r_{k}^5) - 2r_{k}^2\cos(\Theta_0)+2br_{k}^4\sin(\Theta_0)+\mathcal{O}(r_{k}^5)=\nonumber\\
&=& 2(1-\cos(\Theta_0))r_{k}^2+2(a+b\sin(\Theta_0))r_{k}^4+\mathcal{O}(r_{k}^5).
\end{eqnarray}

The term $w_{k}$ is of a form 
$$w_{k}^2=r^2_{k+1}+r^2_{k+q_0}-2r_{k+1}r_{k+q_0}\cos(\varphi_{k+q_0}-\varphi_{k+1}).$$
We use that
$$\varphi_{k+q_0}=\varphi_{k+1}+(q_0-1)\Theta_0 + (q_0-1)br_{k}^2 + \mathcal{O}(r_{k}^4)$$ which we proved in Lemma $\ref{kiteracijaf}$
and
\begin{eqnarray}
\cos(\varphi_{k+q_0}-\varphi_{k+1})&=&\cos((q_0-1)\Theta_0 + (q_0-1)br_{k}^2 + \mathcal{O}(r_{k}^4))=\nonumber\\
&=& \cos((q_0-1)\Theta_0)-(b(q_0-1)r_{k}^2 + \mathcal{O}(r_{k}^4))\sin((q_0-1)\Theta_0)-\nonumber\\
&-&(b(q_0-1)r_{k}^2 + \mathcal{O}(r_{k}^4))^2\frac{\cos((q_0-1)\Theta_0)}{2} + \mathcal{O}(r_{k}^6)=\nonumber\\
&=&\cos((q_0-1)\Theta_0)-b(q_0-1)r_{k}^2\sin((q_0-1)\Theta_0)+\mathcal{O}(r_{k}^4).\nonumber
\end{eqnarray}
It follows
\begin{eqnarray}
w_{k}^2&=&(r_{k}+ar_{k}^3+\mathcal{O}(r_{k}^4))^2+(r_{k}+q_0ar_{k}^3+\mathcal{O}(r_{k}^4))^2-2(r_{k}+ar_{k}^3+\mathcal{O}(r_{k}^4))\nonumber\\
&&(r_{k}+q_0 ar_{k}^3+\mathcal{O}(r_{k}^4))(\cos((q_0-1)\Theta_0)-(q_0-1)br_{k}^2\sin((q_0-1)\Theta_0)+\mathcal{O}(r_{k}^4))=\nonumber\\
&=& 2r_{k}^2 + 2q_0ar_{k}^4+\mathcal{O}(r_{k}^5) - 2r_{k}^2\cos((q_0-1)\Theta_0)+2q_0br_{k}^4\sin((q_0-1)\Theta_0)+\mathcal{O}(r_{k}^5)=\nonumber\\
&=& 2(1-\cos((q_0-1)\Theta_0))r_{k}^2+\mathcal{O}(r_{k}^4).
\end{eqnarray}
We put the behaviours of $\{r_{k}\}$ in the expressions for $y_{k}$, $z_{k}$ and $w_{k}$, and the claim is proven.
The last claim will be showed analogously as in the proof of Lemma 3, only with different coefficient $\sqrt{a^2+b^2}$. $\blacksquare$\\

\textbf{Remark 4.} For the case with $a>0$, we study the inverse system. The lemma is also true for odd $\alpha>3$.\\

By direct calculation we will show in the next theorem the same claim as in the Theorem 3, but for the general rational case and $\alpha=3$.

\begin{theorem}
Let
\begin{eqnarray}
r&\mapsto& r+ar^3+\mathcal{O}(r^4)\nonumber = f(r)\\
\varphi&\mapsto& \varphi + \Theta_0 + br^2 + \mathcal{O}(r^3) =g(\varphi,r)
\end{eqnarray}
be a system with $a\neq 0$, $\Theta_0=\frac{2\pi p}{q}$, for some $p,q\in\mathbf{Z}$ and $e^{i\Theta_0}\neq \sqrt[n]{1}$ for $n=1,2,3,4$. Then there exists $r_1$ small enough such that the discrete spiral $\Gamma(r_0,\varphi_0)$, for $r_0<r_1$ is Minkowski nondegenerate, and $\dim_{B}\Gamma(r_0,\varphi_0)=2/3$.
\end{theorem}
\textbf{Proof.}\\
Without loss of generality, we may assume that $a<0$ and $p=1$. Namely, if $a>0$, then $x_0$ is unstable fixed point. In \cite{laho1}, Lemma 3.1 we showed that the inverse map of $f(r)=r+ar^3+\mathcal{O}(r^4)$ for $a>0$ is of a form $f^{-1}(r)=r-ar^3+\mathcal{O}(r^4)$, so we are back to the case $a<0$.\\
The proof will be done for the case of fixed displacement of angle, that is, $g(\varphi,r)=\varphi+\Theta_0$. In the previous lemma it was shown that in the general case the behaviour of overlapping is the same, and that all other preconditions are valid as well.
The difference in the quotient of the behaviour of $\{z_{k}\}$ can have an influence only on the explicitly calculating Minkowski content, what we will not calculate here.

It follows from the previous lemma the the first overlapping is $z_{k}<2\varepsilon$ if $k>K_0$ for some fixed $K_0$. 
But if we denote the discrete spiral as
$$\Gamma(r_0,\varphi_0)=\cup_{i=0}^{K_0-1}T_{i}(r_{i},\varphi_{i})\cup \Gamma(r_{K_0},\varphi_{K_0}),$$ 
where $T_{i}$ are the points with polar coordinates $(r_{i},\varphi_{i})$, and $\Gamma(r_{K_0},\varphi_{K_0})$ is a discrete spiral with the initial point $(r_{K_0},\varphi_{K_0})$ , that is, $$\Gamma(r_{K_0},\varphi_{K_0})=\cup_{k=K_0}^{\infty}T_{k}(r_{k},\varphi_{k}).$$
Because of the finite stability, it holds that
$$\dim_{B}\Gamma(r_0,\varphi_0)=\dim_{B}\Gamma(r_{K_0},\varphi_{K_0}).$$

Let $\varepsilon_0$ be such that the least number $m_1(\varepsilon_0)$ for which $z_{k}<2\varepsilon_0$ is greater then $K_0$. 
We take fixed $\varepsilon$ such that $\varepsilon<\varepsilon_0$. Now for such $\varepsilon$ we try to find the least natural number $m_1(\varepsilon)$ such that $z_{k}<2\varepsilon$, and we get that $$m_1(\varepsilon)\simeq(q_0a)^{\frac{2}{3}}\varepsilon^{-\frac{2}{3}},$$
that is, 
$$ A_1\varepsilon^{-\frac{2}{3}}\leq m_1(\varepsilon) \leq B_1\varepsilon^{-\frac{2}{3}}.$$
The second overlapping is the overlapping after which we are in the core since we see that the overlappings 
$w_{k}<2\varepsilon$ and $y_{k}<2\varepsilon$ are simultaneous, and we get
$$m_2(\varepsilon)\simeq 2(1-\cos(\Theta_0))\varepsilon^{-2},$$
that is, 
$$A_2\varepsilon^{-2}\leq m_2(\varepsilon) \leq B_2\varepsilon^{-2}.$$

Now we will describe the $\varepsilon$-neighbourhood $A_{\varepsilon}$ of discrete spiral $\Gamma(r_{K_0},\varphi_{K_0})$.
$A_{\varepsilon}$ has three parts:\\
\textbf{1}. $\left|A_1\right|_{\varepsilon}$ - $\varepsilon$-neighbourhood from initial point until first overlapping (\textbf{tail}$_{1}$)\\
\textbf{2}. $\left|A_2\right|_{\varepsilon}$ - $\varepsilon$-neighbourhood Minkowski sausage from first until second overlapping (\textbf{tail}$_{2}$)\\
\textbf{3}. $\left|A_3\right|_{\varepsilon}$ - $\varepsilon$-neighbourhood from second overlapping until the fixed point (\textbf{core})\\

Now we estimate part by part. So,
\begin{equation}
\left|A_1\right|_{\varepsilon}=(m_1(\varepsilon)-K_0)\varepsilon^2\pi\simeq (\varepsilon^{-\frac{2}{3}}-K_0)\varepsilon^2\pi \simeq \varepsilon^{\frac{4}{3}},
\end{equation}
that is,
\begin{equation} \label{kobasmina1}
A_1\pi\varepsilon^{\frac{4}{3}}-K_0\varepsilon^2\pi\leq\left|A_1\right|_{\varepsilon}\leq B_1\pi \varepsilon^{\frac{4}{3}}-K_0\varepsilon^2\pi
\end{equation}
Now the part from first until second overlapping is
\begin{equation} \label{kob2}
\left|A_2\right|_{\varepsilon}=(m_2(\varepsilon)-m_1(\varepsilon))\varepsilon^2\pi - 
\varepsilon^{2}\sum_{k=m_1(\varepsilon)}^{m_2(\varepsilon)-1} (2\arccos(\frac{z_{k}}{2\varepsilon})-\frac{z_{k}}{\varepsilon}\sin(\arccos(\frac{z_{k}}{2\varepsilon})))\\
\end{equation}
We use the equality
$$\sin(\arccos x)=\sqrt{1-x^2}$$
and the following sequences which are valid for $\left|x\right|\leq 1$ of a form
$$\arccos x=\frac{\pi}{2}-a_{0}x-a_{1}x^3 - a_{2}x^5-\ldots - a_{n}x^{2n+1}-\ldots$$
with $a_0=1$ and $a_{n}=\frac{1\cdot3\cdot5\ldots(2n-1)}{2\cdot4\cdot6\ldots(2n)(2n+1)}$ for $n=2,3,\ldots$
and
$$\sqrt{1-x^2}=1+b_1x^2+b_{2}x^4+\ldots+b_{n}x^{2n}+\ldots$$ 
where $b_{i}=-\frac{1\cdot3\cdot5\ldots(2n-3)}{2\cdot4\cdot6\ldots(2n)}$. 

So we simplify the sum
\begin{eqnarray}
&&\sum_{k=m_1(\varepsilon)}^{m_2(\varepsilon)-1} 
(2\arccos(\frac{z_{k}}{2\varepsilon})-\frac{z_{k}}{\varepsilon}\sin(\arccos(\frac{z_{k}}{2\varepsilon})))=\nonumber\\
&=&\sum_{k=m_1(\varepsilon)}^{m_2(\varepsilon)-1} 
[2(\frac{\pi}{2}-\frac{z_{k}}{2\varepsilon}-a_1\frac{z_{k}^3}{8\varepsilon^3}-\ldots-a_{n}\frac{z_{k}^{2n+1}}{(2\varepsilon)^{2n+1}})-\ldots-\\
&-&\frac{z_{k}}{\varepsilon}(1+b_1\frac{z_{k}^2}{4\varepsilon^2}+b_2\frac{z_{k}^4}{2^4\varepsilon^4}+\ldots+b_{n}\frac{z_{k}^{2n}}{2^{2n}\varepsilon^{2n}}+\ldots)]=\nonumber\\
&=&\sum_{k=m_1(\varepsilon)}^{m_2(\varepsilon)-1} (\pi - 2\frac{z_{k}}{\varepsilon}-\frac{a_1+b_1}{4\varepsilon^3}z_{k}^3 -\frac{a_2+b_2}{2^3\varepsilon^5}z_{k}^5-\ldots-\frac{a_{n}+b_{n}}{2^{2n}\varepsilon^{2n+1}}z_{k}^{2n+1}\ldots).\nonumber
\end{eqnarray}
Now we put it in ($\ref{kob2}$) and get
\begin{eqnarray}
\left|A_2\right|_{\varepsilon}=(m_2(\varepsilon)-m_1(\varepsilon))\varepsilon^2\pi-\varepsilon^2\pi\sum_{k=m_1(\varepsilon)}^{m_2(\varepsilon)-1} 1
+2\varepsilon\sum_{k=m_1(\varepsilon)}^{m_2(\varepsilon)-1} z_{k} +\nonumber\\
+\frac{a_1+b_1}{2^2\varepsilon} \sum_{k=m_1(\varepsilon)}^{m_2(\varepsilon)-1} z_{k}^3 + \frac{a_2+b_2}{2^4\varepsilon^3}\sum_{k=m_1(\varepsilon)}^{m_2(\varepsilon)-1}z_{k}^5 +\ldots+\frac{a_{n}+b_{n}}{2^{2n}\varepsilon^{2n-1}}\sum_{k=m_1(\varepsilon)}^{m_2(\varepsilon)-1}z_{k}^{2n+1}+\ldots
\end{eqnarray}
Notice that the first two monomials are abbreviated, so by including the expression $$a_{n}+b_{n}=\frac{2b_{n}}{2n+1}$$ we have
\begin{eqnarray} \label{kobsumea2}
\left|A_2\right|_{\varepsilon}=2\varepsilon\sum_{k=m_1(\varepsilon)}^{m_2(\varepsilon)-1} z_{k} 
+\frac{b_1}{3(2\varepsilon)} \sum_{k=m_1(\varepsilon)}^{m_2(\varepsilon)-1} z_{k}^3 + \frac{b_2}{5(2\varepsilon)^3}\sum_{k=m_1(\varepsilon)}^{m_2(\varepsilon)-1}z_{k}^5+\nonumber\\ +\ldots+\frac{b_{n}}{(2n+1)(2\varepsilon)^{2n-1}}\sum_{k=m_1(\varepsilon)}^{m_2(\varepsilon)-1}z_{k}^{2n+1}+\ldots.
\end{eqnarray}

Now we only need to find a "good" estimations of these sums.

If $A_{z}k^{-\frac{3}{2}}\leq z_{k}\leq B_{z}k^{-\frac{3}{2}}$, then 
$$A_{z}^{2n+1}\sum_{k=a}^{b-1}k^{-(3n+\frac{3}{2})}\leq \sum_{k=a}^{b-1} z_{k}^{2n+1}\leq B_{z}^{2n+1}\sum_{k=a}^{b-1}k^{-(3n+\frac{3}{2})}.$$
We will approximate the sum of powers of $k$ by using the integrals, since it is easy to show that for the decreasing function $f$ it holds
$$\int_{a}^{b} f(x)dx \leq \sum_{k=a}^{b-1} f(k) \leq \int_{a-1}^{b-1} f(x)dx,$$
and we get
$$\int_{k_1}^{k_2} x^{-(3n+\frac{3}{2})}dx = \frac{-1}{3n+\frac{1}{2}}[k_2^{-(3n+\frac{1}{2})}-k_1^{-(3n+\frac{1}{2})}].$$
So, the lower bound for the sum is
$$\sum_{k=m_1(\varepsilon)}^{m_2(\varepsilon)-1} z_{k}^{2n+1}\geq\frac{2A_{z}^{2n+1}}{6n+1}[(\frac{1}{B_1})^{3n+\frac{1}{2}}\varepsilon^{2n+\frac{1}{3}}-(\frac{1}{A_2})^{3n+\frac{1}{2}}\varepsilon^{6n+1}]$$
while the upper bound is
$$ \sum_{k=m_1(\varepsilon)}^{m_2(\varepsilon)-1} z_{k}^{2n+1} \leq \frac{2B_{z}^{2n+1}}{6n+1}[(\frac{1}{A_1})^{3n+\frac{1}{2}}\varepsilon^{2n+\frac{1}{3}}+(3n+\frac{1}{2})(\frac{1}{A_1})^{3n+\frac{3}{2}}\varepsilon^{2n+1}-(\frac{1}{B_2})^{3n+\frac{1}{2}}\varepsilon^{6n+1}+\mathcal{O}(\varepsilon^{2n+\frac{5}{3}})].$$

By using this bounds for power sums of $z_{k}$, we estimate ($\ref{kobsumea2}$) and get
\begin{eqnarray}
\left|A_2\right|_{\varepsilon}\leq \Large[\frac{4B_{z}}{\sqrt{A_1}}+\frac{4A_{z}}{\sqrt{B_1}}\sum_{n=1}^{\infty}\frac{b_{n}}{(2n+1)(6n+1)}(\frac{A_{z}^2}{4B_1^3})^{n} \Large]\varepsilon^{\frac{4}{3}}+ 4B_{z}[(\frac{1}{A_1})^{\frac{3}{2}}-\frac{1}{\sqrt{B_2}}]\varepsilon^2 +\mathcal{O}(\varepsilon^{\frac{8}{3}})\nonumber
\end{eqnarray}
It is easy to show that $\frac{A_{z}^2}{4B_1^3}\leq1$ so the sum 
\begin{equation} \label{koefc1}
C_1=\frac{4B_{z}}{\sqrt{A_1}}+\frac{4A_{z}}{\sqrt{B_1}}\sum_{n=1}^{\infty}\frac{b_{n}}{(2n+1)(6n+1)}\left(\frac{A_{z}}{2\sqrt{B_1^3}}\right)^{2n}
\end{equation}
converges. But we want to show that
$$C_1>0,$$
that is,
\begin{equation} \label{uvjetmin}
\sum_{n=1}^{\infty}\frac{b_{n}}{(2n+1)(6n+1)}(\frac{A_{z}}{2\sqrt{B_1^3}})^{2n}>-\frac{4B_{z}\sqrt{B_1}}{A_{z}\sqrt{A_1}}.
\end{equation}
Recall that $b_{n}<0$. We know that for $x=\frac{A_{z}}{2\sqrt{B_1^3}}<1$
$$\sum_{n=1}^{\infty}\frac{b_{n}}{(2n+1)(6n+1)}x^{2n} >\sum_{n=1}^{\infty}b_{n}x^{2n}=\sqrt{1-x^2}-1>-1.$$
Since $\frac{4B_{z}\sqrt{B_1}}{A_{z}\sqrt{A_1}}>1$, the inequality $(\ref{uvjetmin})$ is true.

Now it follows that 
\begin{equation}
\left|A_2\right|_{\varepsilon} \leq C_1\varepsilon^{\frac{4}{3}} + C_2\varepsilon^{2} +\mathcal{O}(\varepsilon^{\frac{8}{3}})
\end{equation}
where $C_1$ is given by $(\ref{koefc1})$, and $C_2=4B_{z}((\frac{1}{A_1})^{\frac{3}{2}}-\frac{1}{\sqrt{B_2}})$.
Then
\begin{equation} \label{kobmina2}
\frac{\left|A_2\right|_{\varepsilon}}{\varepsilon^{2-s}} \leq C_1\varepsilon^{s-\frac{2}{3}} + C_2\varepsilon^{s} + \mathcal{O}(\varepsilon^{s+\frac{2}{3}}).
\end{equation}

The upper bound for the nucleus is
\begin{equation} \label{kobasmina3}
\left|A_3\right|_{\varepsilon}\leq r_{m_2(\varepsilon)}^2\pi \leq C_3\varepsilon^2.
\end{equation}
So we get the estimation for the Minkowski sausage in the rational case, and it follows from($\ref{kobasmina1}$), ($\ref{kobmina2}$) and ($\ref{kobasmina3}$)
\begin{equation}
\frac{\left|A_1\right|_{\varepsilon}+\left|A_2\right|_{\varepsilon}+\left|A_3\right|_{\varepsilon}}{\varepsilon^{2-s}}\leq (B_1\pi+C_1)\varepsilon^{s-\frac{2}{3}} + (C_2+C_3-\pi)\varepsilon^{s} + \mathcal{O}(\varepsilon^{s+\frac{2}{3}}) \nonumber
\end{equation}
that is,
\begin{equation}
\frac{\left|A\right|_{\varepsilon}}{\varepsilon^{2-s}}\leq (B_1\pi+C_1)\varepsilon^{s-\frac{2}{3}} + (C_2+C_3-\pi)\varepsilon^{s} + \mathcal{O}(\varepsilon^{s+\frac{2}{3}}) \nonumber
\end{equation}
Now we have
$$ \mathcal{M}^{*\frac{2}{3}}(A)\leq B_1\pi+C_1<\infty,$$
and
$$\overline{\rm{dim}}_{B} \Gamma\leq \frac{2}{3}.$$
But for the lower bound we have
\begin{equation}
\frac{\left|A\right|_{\varepsilon}}{\varepsilon^{2-s}}=\frac{\left|A_1\right|_{\varepsilon}+\left|A_2\right|_{\varepsilon}+\left|A_3\right|_{\varepsilon}}{\varepsilon^{2-s}}. \nonumber
\end{equation}

It holds
\begin{equation}
\mathcal{M}^{d}_{*}(A)\geq \mathcal{M}^{d}_{*}(A_1)+\mathcal{M}^{d}_{*}(A_2)+\mathcal{M}^{d}_{*}(A_3).
\end{equation}
Now for $d=\frac{2}{3}$ we have $\mathcal{M}^{\frac{2}{3}}_{*}(A_1)\geq A_1\pi >0$, $\mathcal{M}^{\frac{2}{3}}_{*}(A_2)\geq 0$ and $\mathcal{M}^{\frac{2}{3}}_{*}(A_3)\geq 0$. So we get $\mathcal{M}^{\frac{2}{3}}_{*}(A)>0$, and
$$\underline{\dim}_{B}\Gamma\geq \frac{2}{3}.$$
Finally, it is proven that $\Gamma$ is a Minkowski nondegenerate and $\dim_{B}\Gamma=\frac{2}{3}$. $\blacksquare$\\

\textbf{Remark 5.} It is possible to generalize the claim from the previous proposition for $\alpha=2k+1$ with the change of the nonresonant condition into $e^{i\Theta_0}\neq \sqrt[n]{1}$ for $n=1,2,\ldots,2k+2$, and then $\dim_{B}\Gamma(r_0,\varphi_0)=1-\frac{1}{\alpha}$.

\subsection{The unit-time map of Hopf-Takens bifurcation}

Since the order of overlapping in the irrational case is even more complicated, the direct calculation of the box dimension will be more demanding then for the rational case. In order to avoid it, in this case we are using the connection between the Hopf and Neimark-Sacker bifurcation through the unit-time map. It means that we will prove the box dimension result only for the discrete systems which are unit-time maps of continuous systems with Hopf bifurcation.

At the beginning we have to demonstrate how the unit-time map of Hopf-Takens bifurcation looks like.
We consider the continuous dynamical system
\begin{equation} \label{tok1}
\dot{\mathbf{x}}=\mathbf{F}(\mathbf{x})
\end{equation}
where $\mathbf{x}\in\mathbb{R}^{n}$.
The simplest way to extract the discrete dynamical system from ($\ref{tok1}$) is by using the flow of a system $\phi_{t}(x)$. Namely, we fix $t_0>0$ and we observe the system on $X$ generated by the iteration of  flow $\phi_{t_0}$ (the map with shift $t_0$ along the trajectory of ($\ref{tok1}$)).  Now we chose $t_0=1$ and get the discrete system generated with the unit-time map 
\begin{equation} \label{tok2}
x \mapsto \phi_{1}(x).
\end{equation} 
It is easy to show that the isolated fixed points of ($\ref{tok2}$) corresponds to the isolated singularities of ($\ref{tok1}$). It is also known that the corresponding singularities and fixed points are simultaneously hyperbolic or nonhyperbolic. 

Now, in order to find the unit-time map for the generalized Hopf bifurcation, we start with the normal form for the planar system with the two 
complex conjugated eigenvalues $\pm \omega i$ which has a vector form 
\begin{equation}
\dot{X}=\left(\begin{array}{cc}%
0 & -\omega\\
\omega & 0
\end{array}\right)
\left(\begin{array}{c} %
x\\
y
\end{array}\right) + (x^2+y^2)\left(\begin{array}{cc}%
a & -b\\
b & a
\end{array}\right)\left(\begin{array}{c} %
x\\
y
\end{array}\right),
\end{equation}
i.e. by components,
\begin{eqnarray} \label{hopftak}
\dot{x} &=& -\omega y + (x^2+y^2)(ax-by)\nonumber\\
\dot{y} &=& \omega x + (x^2+y^2) (bx+ay)
\end{eqnarray}
or in polar form
\begin{eqnarray} \label{hopftak33}
\dot{r}&=&ar^3\nonumber\\
\dot{\varphi}&=&\omega + br^2.
\end{eqnarray}
For the system ($\ref{hopftak33}$), the third Lyapunov coefficient is
$$V_3=\frac{9\pi}{\omega} a.$$
So if $a\neq 0$, then $V_3\neq 0$, and by Theorem 5 from $\cite{belg}$ it follows that $\dim_{B}\Gamma=\frac{4}{3}$ where $\Gamma$ is a spiral trajectory near weak focus. Moreover, $\Gamma$ is Minkowski nondgenerate, that is,
\begin{equation}
0<\mathcal{M}_{*}^{d}(\Gamma)\leq\mathcal{M}^{*d}(\Gamma)<\infty
\end{equation}
with $d=\frac{4}{3}$. This result will be used to get the box dimension of Neimark-Sacker bifurcation.
As we know, the unit-time map can be easily obtained by the method of Picard iterations. We get
\begin{eqnarray}
\varphi^{1}(\mathbf{x})=e^{\Lambda}\mathbf{x}+\left\|\mathbf{x}\right\|^2\left(\begin{array}{cc}%
a&-b\\
b&a
\end{array}\right)e^{\Lambda}\mathbf{x} + \mathcal{O}(\left\|\mathbf{x}\right\|^4)
\end{eqnarray}
where
\begin{equation}
e^{\Lambda}=
\left(\begin{array}{cc}%
\cos \omega & -\sin \omega\\
\sin \omega & \cos \omega
\end{array}\right). 
\end{equation}
In the complex form it is
\begin{equation}
f(z) = e^{i\omega}z + \left|z\right|^2(a+ib)e^{i\omega}z + \mathcal{O}(\left|z\right|^4).
\end{equation}
The polar form for $f$ is
\begin{eqnarray} \label{1tokns}
r&\mapsto& r + ar^3 + \mathcal{O}(r^4)\nonumber\\
\varphi & \mapsto & \varphi + \omega + br^2 + \mathcal{O}(r^3).
\end{eqnarray}

Now we can conclude that the normal form map for the Neimark-Sacker bifurcation can be approximated until the 3rd degree by the continuous system which exhibits the Hopf bifurcation. We can also see that all the unit-time maps of the Hopf bifurcation has the form $f$ with the orbits which lie on the trajectories of the system ($\ref{hopftak}$). See Figure 6. By using the property of box dimension regarding the subsets ($A\subseteq B$ $\Rightarrow$ $\dim_{B}A\leq \dim_{B}B$), we get the same upper bound for the box dimension as before, that is, $\frac{4}{3}$. Now we would like to prove that the lower bound is the same, in order to prove that the box dimension is $\frac{4}{3}$.

\begin{center}
\includegraphics[width=4.2cm]{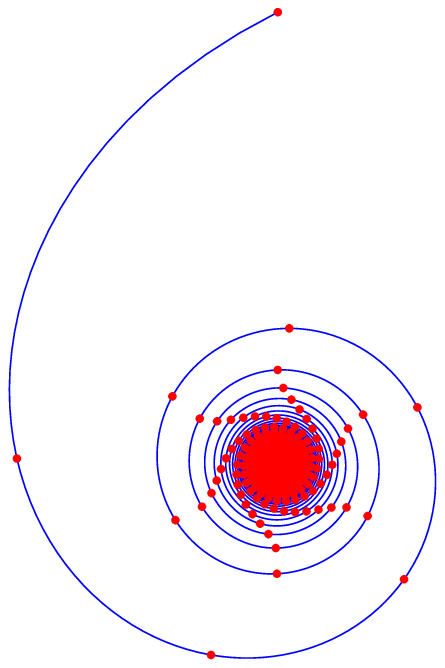}\includegraphics[width=4.2cm]{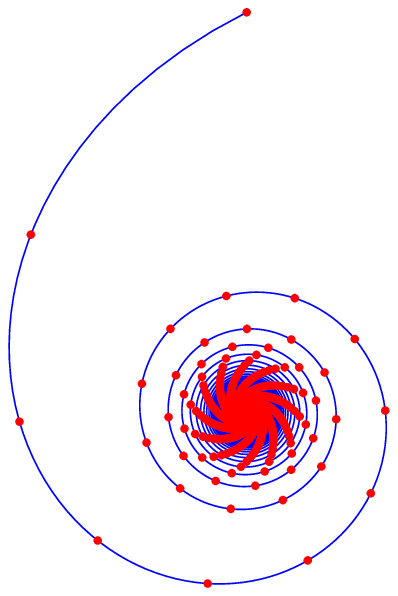}\includegraphics[width=4.2cm]{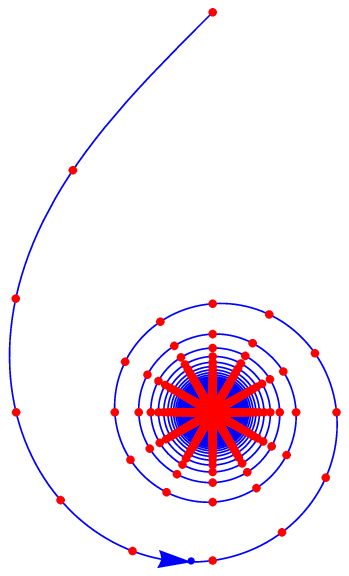}\\
\textit{\textbf{Figure 6a}}\hskip 2cm\textit{\textbf{Figure 6b}}\hskip 2cm  \textit{\textbf{Figure 6c}} \\
\textit{ $a=-1$, $b=1$, $\omega=1$ }\hskip 0.5cm \textit{$a=-1$, $b=1$, $\omega=\frac{\pi}{6}$} \hskip 0.5cm \textit{$a=-1$, $b=0$, $\omega=\frac{\pi}{6}$}
\end{center}

%\textbf{Remark 3.22} 
%It can be shown (for example see ...) that every discrete system of form $(3.93)$  can be approximated until the finite degree by the 1-time map of the continuous system which exhibits the HT bifurcation if the nonresonant condition is satisfied.

Now we observe the generalized planar continuous system 
\begin{eqnarray}
\dot{x}&=&-\omega y + (x^2+y^2)^{k}(ax-by)\nonumber\\
\dot{y}&=&\omega x + (x^2+y^2)^{k}(bx+ay)
\end{eqnarray}
which in polar coordinates is 
\begin{eqnarray}
\dot{r}&=&ar^{2k+1}\nonumber\\
\dot{\varphi}&=&\omega+br^{2k}.
\end{eqnarray}

\begin{lemma}
If we have the planar system
\begin{eqnarray}
\dot{r}&=&a r^{2k+1}\nonumber\\
\dot{\varphi}&=&\omega+b r^{2k},
\end{eqnarray}
then the corresponding unit-time map is of a form
\begin{eqnarray}
r&\mapsto& r+ar^{2k+1}+\mathcal{O}(r^{2k+2}) \nonumber\\
\varphi&\mapsto& \varphi + \omega +br^{2k} + \mathcal{O}(r^{2k+1}).
\end{eqnarray}
\end{lemma}

\textbf{Remark 6.} The lemma can be easily proven using the Picar iteration method.

\textbf{Remark 7.} By adding the terms of higher order to the continuous system, the form of the unit-time map doesn't change.

The following lemma is a generalisation of the box dimension result from \cite{zuzu} (Theorem 9(b)) and it can be shown analogously. 

\begin{lemma}
Let $\Gamma$ be a spiral trajectory near the origin of planar system
\begin{eqnarray} \label{le1}
\dot{r}&=&a r^{2k+1}\nonumber\\
\dot{\varphi}&=&\omega+b r^{2k}.
\end{eqnarray}
Then $\Gamma$ is a Minkowski nondegenerate and we have
$$\dim_{B}\Gamma=2(1-\frac{1}{2k+1}).$$
\end{lemma}
\textbf{Proof.}\\
This lemma will be proven by using Theorem 7 from \cite{zuzu}. It means that we have to show that $r=f(\varphi)$ is a radially decreasing function, and that it satisfies the following conditions: $r=f(\varphi)\simeq \varphi^{-\frac{1}{2k}}$, $f(\varphi)-f(\varphi +2\pi)\simeq  \varphi^{-\frac{1}{2k}-1}$ and $\left|f'(\varphi)\right|\simeq \varphi^{-\frac{1}{2k}-1} $ and $\left|f''(\varphi)\right|\leq M_3\varphi^{-\frac{1}{2k}}$.
From the assumption $a<0$, it follows that $r=f(\varphi)$ is a radially decreasing function, that is, the origin is a stable focus. 

Now we seek for the solution of a system ($\ref{le1}$), and get
$$\frac{\dot{r}}{\dot{\varphi}}=\frac{ar^2k+1}{\omega + br^{2k}}$$
or
$$\frac{\omega + br^{2k}}{ar^{2k+1}}dr=d\varphi.$$
By integrating, it follows
$$\frac{-\omega}{a(2k)}r^{-2k} + \frac{b}{a}\ln(r)=\varphi+C,$$
ie.
$$\varphi=\Phi(r)=\frac{-\omega}{a(2k)}r^{-2k} + \frac{b}{a}\ln(r) - C.$$

It is easy to show that $\varphi=\Phi(r)\simeq r^{-2k}$, and then we have $r=f(\varphi)\simeq \varphi^{-\frac{1}{2k}}$.
Then $\left|f'(\varphi)\right|\simeq \varphi^{-\frac{1}{2k}-1}$. 
By Mean Value Theorem we get
$f(\varphi)-f(\varphi +2\pi)\simeq  \varphi^{-\frac{1}{2k}-1}$.
We know
$$\Phi''(r)=-\frac{f''(\Phi(r))\Phi'(r)}{(f'(\Phi(r)))^{2}},$$
so it follows the last condition on the derivative
$$\left|f''(\varphi)\right|=\left|f'(\varphi)^2 \Phi''(r) \Phi'(r)^{-1}\right|\leq c \varphi^{-\frac{1}{k}-2}\varphi^{-1-\frac{1}{2k}}\varphi^{1+\frac{1}{k}}\leq M_3\varphi^{-\frac{1}{2k}}.$$
Now we have $\dim_{B}\Gamma=2(1-\frac{1}{2k+1})$.
$\blacksquare$\\

Until now, we observed the truncated normal form for the Hopf bifurcation at the bifurcation value $\mu=0$ because then the origin is a nonhyperbolic singularity and the box dimension is nontrivial. But, if we consider the whole one-parameter family of planar vector fields which exhibits the Hopf bifurcation at $\mu=0$, it is of a form
\begin{eqnarray} \label{opethopf}
\dot{x}&=& d\mu x - (\omega +c\mu)y + (ax-by)(x^2+y^2)\nonumber\\
\dot{y}&=& (\omega+c\mu)x + d\mu y + (bx+ay)(x^2+y^2)
\end{eqnarray}
or in polar coordinates
\begin{eqnarray}
\dot{r}&=&d\mu r + ar^3\nonumber\\
\dot{\varphi}&=& \omega +c\mu + br^2
\end{eqnarray}
with the transversality condition $d\neq 0$ and nondegeneracy condition $a\neq 0$.

Analogously as before, we get the unit-time map
\begin{eqnarray} \label{opet1tok}
r&\mapsto& r+d\mu r +ar^3+ \mathcal{O}(\mu^2r,\mu r^2,r^4)\nonumber\\
\varphi&\mapsto& \varphi + \omega + c\mu+ br^2 + \mathcal{O}(\mu^2,\mu r^2,r^4).
\end{eqnarray}
It is known that the transversality condition for the Neimark-Sacker bifurcation is also $d\neq 0$, and the nondegeneracy condition is also $a\neq 0$.
Since the nonresonant condition for the Neimark-Sacker bifurcation ensures that the normal form is given by $(\ref{opet1tok})$, we don't need this condition since every unit-time map is of a given form.
So, in the same time the system ($\ref{opethopf}$) experiences the Hopf bifurcation for $\mu=0$ and the corresponding unit-time map experiences the Neimark-Sacker bifurcation for $\mu=0$. At Figure 7 we can see the Hopf bifurcation and Neimark-Sacker bifurcation with $\omega=1$, $a=d=1$ and $b=c=0$. 

\begin{center}
\includegraphics[width=3.5cm]{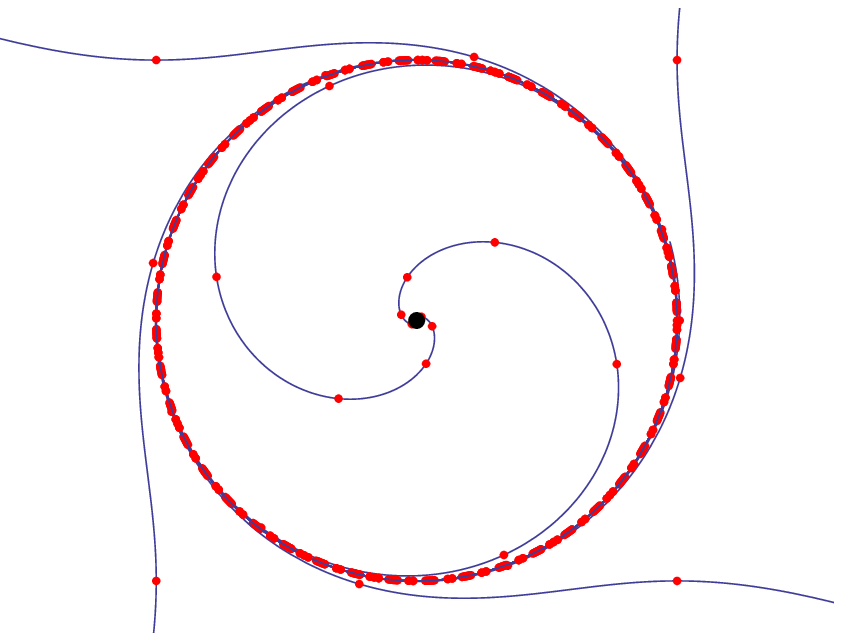}\hskip 1cm\includegraphics[width=3.5cm]{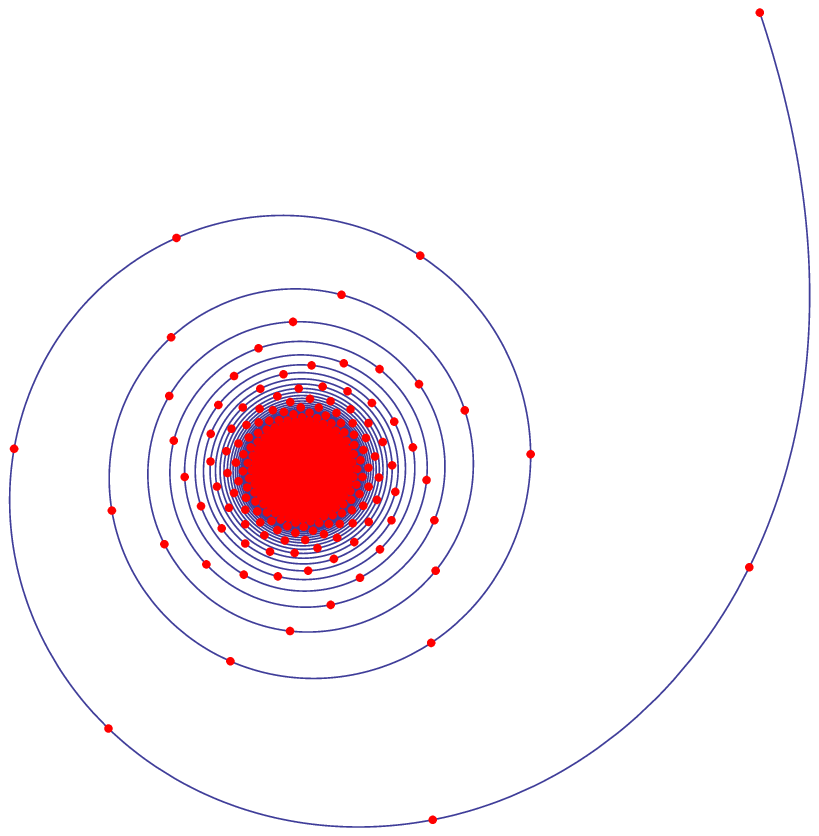}\hskip 1cm\includegraphics[width=3cm]{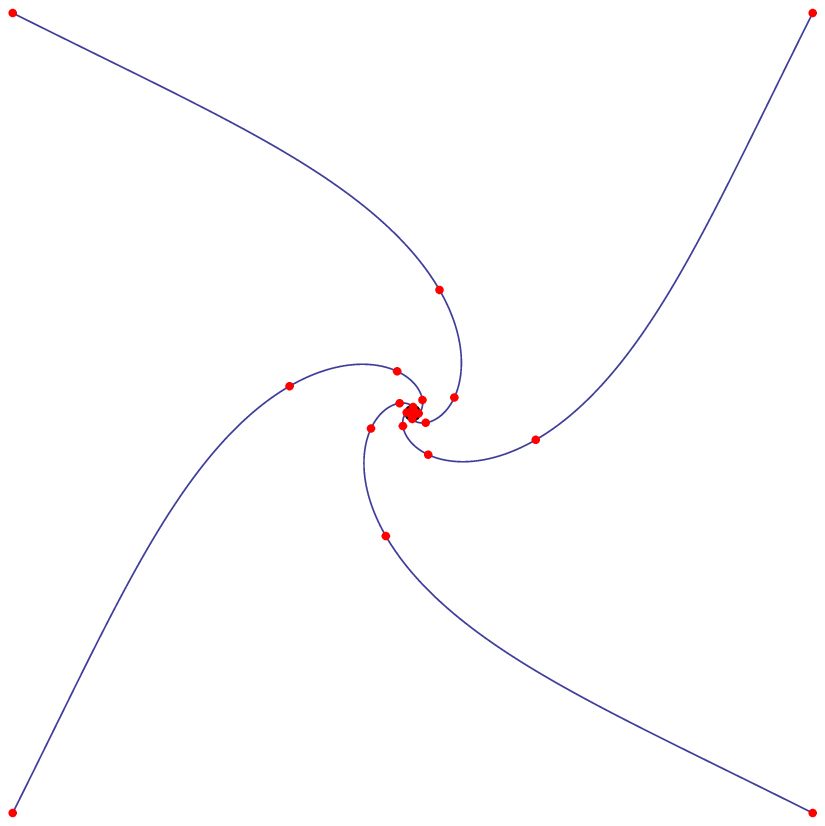}\\
\textit{\textbf{Figure 7a} $\mu>0$} \hskip 2cm \textit{\textbf{Figure 7b} $\mu=0$ }\hskip 2cm \textit{\textbf{Figure 7c} $\mu>0$}
\end{center}

\subsection{Irrational case}

In this chapter we will calculate the box dimension of the orbit of discrete system which exhibits the Neimark-Sacker bifurcation in irrational case. The normal form at the bifurcation value is
\begin{eqnarray} \label{1tokns}
r&\mapsto& r + ar^3 + \mathcal{O}(r^5)\nonumber\\
\varphi & \mapsto & \varphi + \Theta_0 + br^2 + \mathcal{O}(r^4),
\end{eqnarray}
where $\Theta_0=2\pi\beta$, $\beta\in\mathbb{R}\backslash \mathbb{Q}$ (irrational case). 

Before we begin with the ordering of overlappings, notice that the angle mapping depends on the radius. So we observe the expressions for $y_{k}$, $z_{k}$ and $w_{k}$. Since $\Theta_0$ is irrational, then there exists the lowest natural number $q_0$ such that $q_0\Theta_0>2\pi$, that is $q_0=\left\lceil \frac{2\pi}{\Theta_0}\right\rceil$. 
Let us look at the behaviour of the sequences $\{y_{k}\}$, $\{z_{k}\}$ i $\{w_{k}\}$ defined in Subsection 3.2.
\begin{lemma}
Let
$$r_{k+1}=r_{k}+ar_{k}^3+\mathcal{O}(r_{k}^5)$$
$$\varphi_{k+1}=\varphi_{k} + \Theta_0 + br_{k}^2 +\mathcal{O}(r_{k}^4)$$
be a given system with  $a\neq 0$.
Then for $\Theta_0=2\pi\beta$, $\beta\in\mathbb{R}\backslash\mathbb{Q}$ where $q_0=\left\lceil \frac{2\pi}{\Theta_0}\right\rceil$ we have
\begin{eqnarray}
y^2_{k} &=& r_{k}^2+r_{k+1}^2 - 2r_{k}r_{k+1}\cos(\Theta_0)+\mathcal{O}(r_{k}^4) = 2(1-\cos(\Theta_0)) r_{k}^2 + \mathcal{O}(r_{k}^4)\nonumber\\
z^2_{k} &=& r_{k}^2+r_{k+q_0}^2 - 2r_{k}r_{k+q_0}\cos(q_0\Theta_0)+\mathcal{O}(r_{k}^4) = 2(1-\cos(q_0\Theta_0))r_{k}^2 + \mathcal{O}(r_{k}^4)\nonumber\\
w^2_{k} &=& r_{k+1}^2+r_{k+q_0}^2 - 2r_{k+1}r_{k+q_0}\cos((q_0-1)\Theta_0)+\mathcal{O}(r_{k}^4) =\nonumber\\
&=& 2(1-\cos(q_0-1)\Theta_0)r_{k}^2  +\mathcal{O}(r_{k}^4).\nonumber
\end{eqnarray}
\end{lemma}

\textbf{Proof.}\\
The behaviuor of $\{y_{k}\}$ can be showed as in Lemma 3.
We consider $$z_{k}^2=r_{k}^2 + r^2_{k+q_0}-2r_{k}r_{k+q_0}\cos(\varphi_{k+q_0}-\varphi_{k}).$$ 
Using the expressions for $r_{k+q_0}$ and $\varphi_{k+q_0}$, we get 
\begin{eqnarray}
z_{k}^2&=&r_{k}^2+(r_{k}+q_0ar_{k}^3+\mathcal{O}(r_{k}^4))^2-\nonumber\\
&-&2r_{k}(r_{k}+q_0 ar_{k}^3+\mathcal{O}(r_{k}^4))(\cos(q_0\Theta_0)-q_0br_{k}^2\sin(q_0\Theta_0)+\mathcal{O}(r_{k}^4))=\nonumber\\
&=& 2r_{k}^2 + 2q_0ar_{k}^4+\mathcal{O}(r_{k}^5) - 2r_{k}^2\cos(q_0\Theta_0)+2q_0br_{k}^4\sin(q_0\Theta_0)+\mathcal{O}(r_{k}^5)=\nonumber\\
&=& 2(1-\cos(q_0\Theta_0))r_{k}^2+2q_0(a+b\sin(q_0\Theta_0))r_{k}^4+\mathcal{O}(r_{k}^5).
\end{eqnarray}
Analogously we get
\begin{eqnarray}
w_{k}^2&=&(r_{k}+ar_{k}^3+\mathcal{O}(r_{k}^4))^2+(r_{k}+q_0ar_{k}^3+\mathcal{O}(r_{k}^4))^2-2(r_{k}+ar_{k}^3+\mathcal{O}(r_{k}^4))\nonumber\\
&&(r_{k}+q_0 ar_{k}^3+\mathcal{O}(r_{k}^4))(\cos((q_0-1)\Theta_0)-(q_0-1)br_{k}^2\sin((q_0-1)\Theta_0)+\mathcal{O}(r_{k}^4))=\nonumber\\
&=& 2r_{k}^2 + 2q_0ar_{k}^4+\mathcal{O}(r_{k}^5) - 2r_{k}^2\cos((q_0-1)\Theta_0)+\nonumber\\
&+&2(q_0-1)br_{k}^4\sin((q_0-1)\Theta_0)+\mathcal{O}(r_{k}^5)=\nonumber\\
&=& 2(1-\cos((q_0-1)\Theta_0))r_{k}^2+\mathcal{O}(r_{k}^4).
\end{eqnarray}
$\blacksquare$\\

Before we start to calculate the $\varepsilon$-neighbourhood $A_{\varepsilon}$ of the orbit of the system ($\ref{1tokns}$), we are interested in the order of overlappings of neighboring terms of sequences $y_{k}$, $w_{k}$ and $z_{k}$ since all of them behave like $k^{-\frac{1}{2}}$ but with different coefficients. 
If the first overlapping is between the adjoining points ($y_{k}<\varepsilon$), then the lower bound for $A_{\varepsilon}$ of that overlapping is in fact the $\varepsilon^{*}$-neighbourhood of the corresponding spiral trajectory with the radius $\varepsilon^{*}>\frac{\sqrt{3}-1}{2}\varepsilon$. That will be the part of the tail for the $\varepsilon^{*}$-neighbourhood $\left|A_{\varepsilon^{*}}\right|$ for spiral trajectory. Furthermore, if for the spiral trajectory with $\varepsilon^{*}$ the entering in the nucleus is after the overlapping $y_{k}<\varepsilon$, then that nucleus with $\varepsilon^{*}$ is the lower bound for the nucleus with $\varepsilon$.

Therefore, if the first overlapping is the overlapping of adjoining points, then we have the result as we will see in the proof of the theorem . Now we want to show that this is always true.

\begin{lemma}
Let be given a system ($\ref{1tokns}$) with $\Theta_0=2\pi\beta$, $\beta\in\mathbb{R}\backslash\mathbb{Q}$. Then there exist $\Theta_1$ small enough and $K_0$ big enough such that for all $\Theta_0<\Theta_1$ and $k>K_0$, it holds $y_{k}<z_{k}$.
\end{lemma}

\textbf{Proof.}\\
By considering the behaviour of the sequence
$$y_{k}=d(A_{k},A_{k+1})=\sqrt{r^2_{k}+r^2_{k+1}-2r_{k}r_{k+1}\cos(\Theta_0)+\mathcal{O}(r_{k}^4)},$$
we see that
$$y^2_{k} = (r_{k}-r_{k+1})^2 + 2r_{k}r_{k+1}(1-\cos(\Theta_0))+\mathcal{O}(r_{k}^4)>(r_{k}-r_{k+1})^2 + \frac{d_1}{k^2}.$$
Notice that if the angle $\Theta_0$ is smaller, then $y_{k}$ is closer to $(r_{k}-r_{k+1})+\frac{d_1}{k^2}$, ie. when $\Theta_0 \rightarrow 0$, then $y_{k}\rightarrow (r_{k}-r_{k+1})$. By considering the behaviour of 
$$z_{k}=d(A_{k},A_{k+q_0})=\sqrt{r^2_{k}+r^2_{k+q_0}-2r_{k}r_{k+q_0}\cos(q_0\Theta_0)+\mathcal{O}(r_{k}^4)}$$
we notice that
$$z^2_{k}=(r_{k}-r_{k+q_0})^2+2r_{k}r_{k+q_0}(1-\cos(q_0\Theta_0))+\mathcal{O}(r_{k}^4)>(r_{k}-r_{k+q_0})^2 + \frac{d_2}{k^2}.$$
We see that if $\Theta_0$ is smaller, then $q_0\Theta_0(\rm{mod} 2\pi)$ is smaller, so $z_{k}$ is closer to $r_{k}-r_{k+q_0}+\frac{d_2}{k^2}$. But also when $\Theta_0\rightarrow 0$, then $q_0\rightarrow \infty$ and $r_{k}-r_{k+q_0}$ increases from $r_{k}-r_{k+1}+\frac{d_2}{k^2}$ to $r_k+\frac{d_2}{k^2}$. 
It means that we get: if $\Theta_0\rightarrow 0$, then $y_{k}$ decreases to $r_{k}-r_{k+1}+\frac{d_1}{k^2}$, and $z_{k}$ increases from $r_{k}-r_{k+1}+\frac{d_2}{k^2}$ to $r_{k}+\frac{d_2}{k^2}$. 
Hence, we can chose $\Theta_0$ small enough such that $y_{k}<z_{k}$, for every $k$ big enough. 
Let $\Theta_1$ be such an angle, ie. for all $\Theta_0<\Theta_1$ we have $y_{k}<z_{k}$, for $k$ big enough. $\blacksquare$\\

Theorem 2 showed that
$$\frac{2}{3}\leq \underline{\dim}_{B}\Gamma\leq\overline{\dim}_{B}\Gamma\leq\frac{4}{3}.$$
It is already proven that for $\Theta_0=2\pi\beta$, $\beta\in\mathbb{Q}$ the box dimension is equal to the lower bound. 
In the next theorem we will prove that in the irrational case the box dimension is equal to the upper bound.

\begin{theorem} \textbf{Box dimension for Neimark-Sacker bifurcation - irrational case}\\
Let we have a discrete planar system
\begin{eqnarray} \label{nsirac}
r&\mapsto& r + ar^3 + \mathcal{O}(r^5) = f(r)\nonumber\\
\varphi & \mapsto & \varphi + \Theta_0 + br^2 + \mathcal{O}(r^4) = g(r,\varphi).
\end{eqnarray}
which is the unit-time map of the continuous planar system
\begin{eqnarray} \label{hopftakpol}
\dot{r}&=&ar^3 + \mathcal{O}(r^5)\nonumber\\
\dot{\varphi}&=&\Theta_0 + br^2 +\mathcal{O}(r^4),
\end{eqnarray}
and let $\Theta_0=2\pi\beta$, with $\beta$ irrational. Then there exists $r_0>0$ such that the sequence $S(r_1,\varphi_1)=\{(r_{k},\varphi_{k})\}$ defined by $r_{k+1}=f(r_{k})$ and $\varphi_{k+1}=g(r_{k},\varphi_{k})$ 
where $r_1<r_0$, is a Minkowski nondegenerate , and we have
$$\dim_{B}S(r_1,\varphi_1)=\frac{4}{3}.$$
\end{theorem}

\textbf{Proof.}\\

First we should check the ordering of overlappings. If $\omega<\Theta_1$, then it follows from Lemma 8 that $y_{k}<z_{k}$ for $k$ big enough. So in that case the first overlapping is for the sequence $y_{k}$. But, what if the condition $\omega>\Theta_1$ isn't satisfied. Then we would like to find some other system which is also the $T_0$-time map of ($\ref{hopftakpol}$), but with $\omega_{*}$ less then $\Theta_1$. And, of course, we would like that system has the same box dimension as the initial system.
We consider the system
\begin{eqnarray}
r&\mapsto&f_1(r)\nonumber\\
\varphi &\mapsto& g_1(r,\varphi)
\end{eqnarray} 
where $f_1$ and $g_1$ are such that $f_1^{q}=f$ i $g_1^{q}=g$ for the lowest natural number $q$ such that $\frac{\Theta_0}{q}<\Theta_1$.
It is easy to see that the previous system is of a form 
\begin{eqnarray} \label{qnsiracq}
r&\mapsto& r + \frac{a}{q}r^3 + \mathcal{O}(r^5) = f(r)\nonumber\\
\varphi & \mapsto & \varphi + \frac{\Theta_0}{q} + \frac{b}{q}r^2 + \mathcal{O}(r^4) = g(r,\varphi).
\end{eqnarray}
Furthermore, notice that since the system ($\ref{nsirac}$) is the unit-time map of the system ($\ref{hopftakpol}$), then in fact the system ($\ref{qnsiracq}$) is a $\frac{1}{q}$-time map of the same continuous system. It remains to show that the box dimensions of the orbits of systems ($\ref{nsirac}$) and ($\ref{qnsiracq}$) are the same. Namely, if $\Gamma^{*}$ is a discrete spiral of ($\ref{qnsiracq}$) with the initial point $(r_0,\varphi_0)$, then
$$\Gamma^{*}(r_0,\varphi_0)=\cup_{i=1}^{q}\Gamma_{i}(r_0^{i},\varphi_0^{i}),$$
where $\Gamma_i$ are the discrete spirals of ($\ref{nsirac}$) with the initial point $$(r_0^{i},\varphi_0^{i})=(f_1^{i-1}(r_0),g_1^{i-1}(r_0,\varphi_0)).$$ Now, it follows from ($\ref{konstab}$)
$$\dim_{B} \Gamma^{*}(r_0,\varphi_0)=\dim_{B}\Gamma(r_0,\varphi_0),$$
where $\Gamma$ is a discrete spiral of ($\ref{nsirac}$) with the initial point $(r_0,\varphi_0)$.\\

Hence, without loss of generality we may assume that $\omega<\Theta_1$, that is, that the first overlapping is $y_{k}<2\varepsilon$. Furthermore, since we would like to find the lower bound, we can take the inequality $y_{k}<\varepsilon$ instead of $y_{k}<2\varepsilon$. The difference in the estimation of $\varepsilon$-neighbourhood will not influence the lower bound. Then we get the lower bound for the corresponding spiral trajectory by using the "annular" spiral. See Figure 8.
If $y_{k}<\varepsilon$, then the radius of "annulus" is $\varepsilon^{*}$ such that $\varepsilon^{*}>\frac{\sqrt{3}-1}{2}\varepsilon$.\\
We denote by $A=\Gamma(r_0,\varphi_0)$ discrete spiral, and with $\Gamma$ spiral trajectory of the system ($\ref{hopftakpol}$).

\begin{center}
\includegraphics[width=6cm]{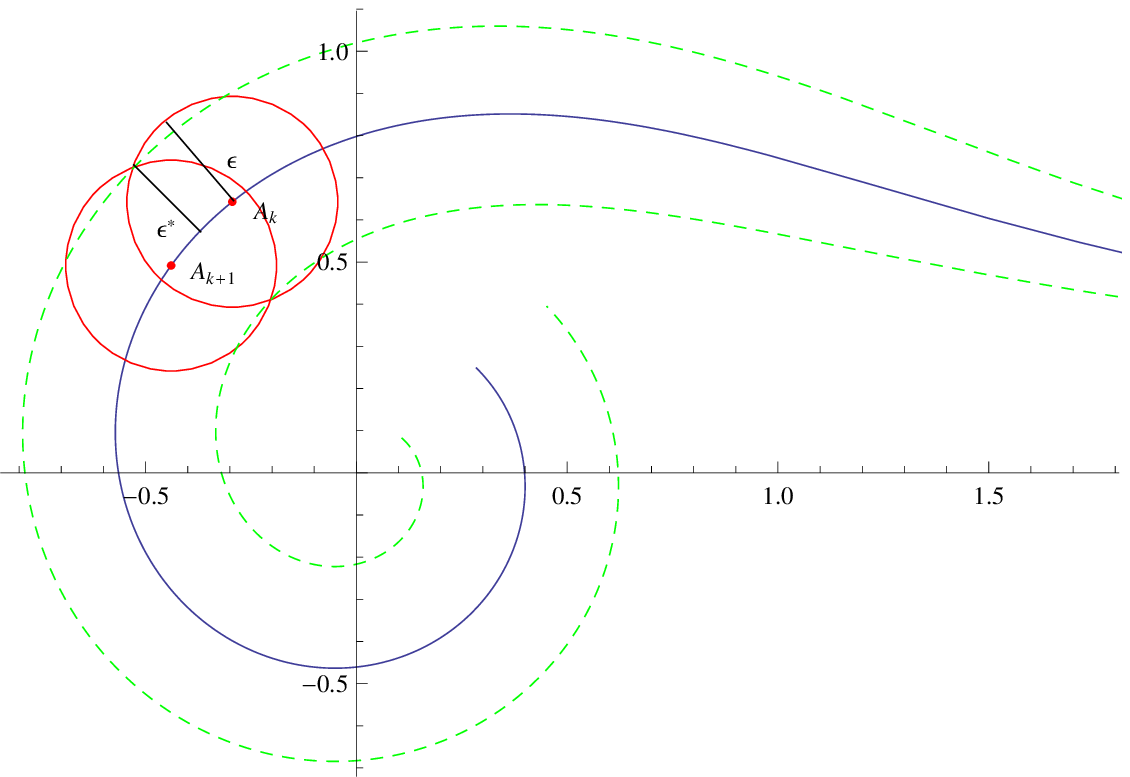}\\
\textit{\textbf{Figure 8}} \\
\end{center}

We divide $\varepsilon$-neighbourhood $A_{\varepsilon}$ into several pieces to simplify the estimation process.\\

\textbf{\textit{PART 1.}} $A_{1,\varepsilon}$ = $\varepsilon$-neighbourhood from $(r_0,\varphi_0)$ until the first overlapping
\begin{eqnarray}
\left|A_{1,\varepsilon}\right|&=&\varepsilon^2\pi(m_1(\varepsilon)-1)\geq \varepsilon^2\pi(c_{1d}\varepsilon^{-\frac{2}{3}}-1)=\pi(c_{1d}\varepsilon^{\frac{4}{3}}-\varepsilon^2)\\
\frac{\left|A_{1,\varepsilon}\right|}{\varepsilon^{2-s}}&\geq& \pi(c_{1d}\varepsilon^{s-\frac{2}{3}}-\varepsilon^{s})
\end{eqnarray}
where $m_1(\varepsilon)\geq c_{1d}\varepsilon^{-\frac{2}{3}}$ since $y_{k}\geq (r_{k}-r_{k+1})\simeq k^{-\frac{3}{2}}$.\\

\textbf{\textit{PART 2.}} $A_{2,\varepsilon}$ = $\varepsilon$-neighbourhood from the first overlapping until the end consist of two parts: the tail $\left|A_{2,\varepsilon}\right|$ which contains all possible overlappings and the nucleus $\left|A_{3,\varepsilon}\right|$. This part will be estimated by the $\varepsilon^{*}$-neighbourhood with $\varepsilon^{*}$ from the corresponding spiral trajectory of continuous systems near focus which also contains two different parts:
the tail  $\left|\Gamma_{r,\varepsilon^{*}}\right|$
from $\varphi_{m_1(\varepsilon)}$ do $\varphi_2(\varepsilon^{*})$ (this is the angle after which we are in the nucleus) and the nucleus $\left|\Gamma_{j,\varepsilon_{*}}\right|$.
Hence, it holds 
\begin{eqnarray}
\left|A_{2,\varepsilon}\right|+\left|A_{3,\varepsilon}\right|>\left|\Gamma_{r,\varepsilon^{*}}\right|+\left|\Gamma_{j,\varepsilon_{*}}\right|
\end{eqnarray}
Notice that the right side of $\varepsilon$-neighbourhood is continuous spiral $\Gamma$ for $\varphi\in[\varphi_{m_1(\varepsilon)},\infty )$. 
Now we use that $\Gamma$ is a Minkowski nondegenerate and for $d=\frac{4}{3}$ we have
$$\mathcal{M}^{d}_{*}(\Gamma)=\liminf_{\varepsilon\rightarrow 0} \frac{\left|\Gamma_{\varepsilon}\right|}{\varepsilon^{2-d}}>0.$$
It follows
\begin{equation}
\frac{\left|A_{2,\varepsilon}\right|+\left|A_{3,\varepsilon}\right|}{\varepsilon^{2-s}}\geq \frac{\left|\Gamma_{r,\varepsilon^{*}}\right|+\left|\Gamma_{j,\varepsilon^{*}}\right|}{\varepsilon^{2-s}}
\end{equation}
Now we include $s=d$ to get 
\begin{equation}
\frac{\left|A_{2,\varepsilon}\right|+\left|A_{3,\varepsilon}\right|}{\varepsilon^{2-d}}\geq \frac{\left|\Gamma_{r,\varepsilon^{*}}\right|+\left|\Gamma_{j,\varepsilon^{*}}\right|}{\varepsilon^{2-d}}.
\end{equation}
It follows from $\varepsilon^{*}>\frac{\sqrt{3}-1}{2}$ that $\varepsilon^{d-2}>\left(\frac{\sqrt{3}-1}{2}\right)^{2-d}(\varepsilon^{*})^{d-2}$, so we get
\begin{equation}
\frac{\left|A_{2,\varepsilon}\right|+\left|A_{3,\varepsilon}\right|}{\varepsilon^{2-d}}\geq (\frac{\sqrt{3}-1}{2})^{2-d}\frac{\left|\Gamma_{r,\varepsilon^{*}}\right|+\left|\Gamma_{j,\varepsilon^{*}}\right|}{(\varepsilon^{*})^{2-d}}.
\end{equation}

Now, taking into consideration the parts 1. and 2. we have
\begin{equation}
\mathcal{M}_{*}^{d}(A)=\liminf_{\varepsilon\rightarrow 0}\frac{\left|A_{1,\varepsilon}\right|+\left|A_{2,\varepsilon}\right|+\left|A_{3,\varepsilon}\right|}{\varepsilon^{2-d}}
\geq 0+\left(\frac{\sqrt{3}-1}{2}\right)^{2-d}\mathcal{M}_{*}^{d}(\Gamma)>0.
\end{equation}
Therefore, we have the lower inequality for the $\varepsilon$-neighbourhood.\\

On the other hand, it is quite obvious that the $\varepsilon$-neighbourhood for the discrete spiral $A_{\varepsilon}$ is contained in the  $\varepsilon$-neighbourhood of the corresponding continuous spiral $\Gamma_{c,\varepsilon}$ ie.
\begin{equation}
\left|A_{1,\varepsilon}\right|+\left|A_{2,\varepsilon}\right|+\left|A_{3,\varepsilon}\right|\leq\left|\Gamma_{ct,\varepsilon}\right|+\left|\Gamma_{cn,\varepsilon}\right|=\left|\Gamma_{c,\varepsilon}\right|
\end{equation} 
where $\Gamma_{ct,\varepsilon}$ is a tail of $\varepsilon$-neighbourhood of continuous spiral $\Gamma$, and $\Gamma_{cn,\varepsilon}$ is a corresponding nucleus.
So we have
\begin{equation}
\mathcal{M}^{*d}(A)=\limsup_{\varepsilon\rightarrow 0} \frac{\left|A_1\right|_{\varepsilon}+\left|A_{21}\right|_{\varepsilon}+
\left|A_{22}\right|_{\varepsilon}}{\varepsilon^{2-d}}\leq\mathcal{M}^{*d}(\Gamma)<\infty. 
\end{equation}

Hence, we have
$$0<\mathcal{M}_{*}^{d}(A)\leq \mathcal{M}^{*d}(A)<\infty$$
and it follows that 
$$\dim_{B}(A)=d.\,\,\,\,\blacksquare$$

\textbf{Remark 6.} Using the Lemma 5 and Lemma 6, it is possible to prove the analogous result for the generalized Neimark-Sacker bifurcation, that is, Chenciner bifurcation (see \cite{ch1}, \cite{ch2}, \cite{kuz},). If we apply the obtained result for the box dimension of the Neimark-Sacker bifurcation in irrational case to Chenciner bifurcation, we get that the box dimension of the orbit near the origin at the bifurcation point $(\beta_1,\beta_2)=(0,0)$ is $\dim_{B}S(r_1,\varphi_1)=2(1-\frac{4}{5})=\frac{8}{5}$, as we can see at Figure 9.

\begin{center}
\includegraphics[width=6cm]{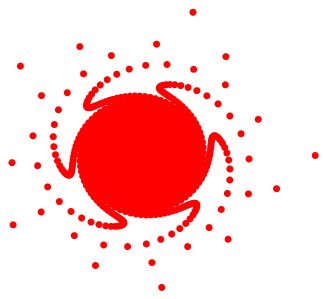}\includegraphics[width=6cm]{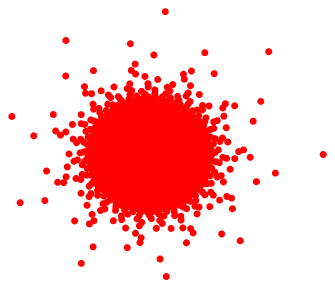}\\
\textit{\textbf{Figure 9a} $f(r)=r-r^5$} \hskip 3cm \textit{\textbf{Figure 9b} $f(r)=r-r^5$}\\
\textit{$g(\varphi,r)=1+r^2$}\hskip 3cm \textit{$g(\varphi,r)=0.5+r^2$}\\
\end{center}

\textbf{Remark 9.} Regarding the Neimark-Sacker bifurcation, we can conclude that the difference in the box dimension is connected to the different order of overlappings. Namely, in the rational case it is $z_{k}<y_{k}$ for all $k$ big enough, while in the irrational case the ordering is reverse. 
We can also see that the Neimark-Sacker bifurcation has a nondegeneracy condition of order three. So, like in all previously studied bifurcation (see \cite{laho1},\cite{zuzu}), the box dimension is connected with the order of nondegeneracy. But in the case of this bifurcation, it is also connected to the rationality of displacement angle. Moreover, notice that the discrete spiral in rational case has the same box dimension as Poincar$\acute{e}$ map of continuous spiral near weak focus (see \cite{belg}), while in the irrational case we get the box dimension of continuous spiral. It is also interesting to notice that the box dimension of bifurcated object differs from $0$ (in the case of periodic orbit on the invariant circle) to $1$ (in the case of dense orbit on the invariant circle).

\textbf{Remark 10.} It will be interesting to further explore Neimark-Sacker bifurcation of limit cycle and the bifurcation of the systems in higher dimension with Neimark-Sacker bifurcation such as fold-Neimark-Sacker (fold-NS, see \cite{bsv} ), flip-Neimark-Sacker (flip-NS, see \cite{kuz2}) and double Neimark-Sacker (NS-NS, see \cite{kuz2}). \\

\textbf{Remark 11.} This connection via unit-time map between the discrete and continuous dynamical systems can be applied to fractal analysis of nilpotent singularities of planar vector fields. See \cite{hozu}.\\

\textbf{Acknowledgements.}\\
 I would like to thank Prof. V. \v Zupanovi\' c for helpful comments and suggestions.

\textit{E-mail address}: lana.horvat@fer.hr

University of Zagreb,\\
Faculty of Electrical Engineering and Computing,\\ 
Department of Applied Mathematics, \\
Unska 3, 10000 Zagreb, Croatia\\


\begin{thebibliography}{99}

\bibitem{ap} D.K. Arrowsmith, C.M. Place, \textit{An Introduction to dynamical systems}, Cambridge University Press, Cambridge, 1990

\bibitem{bave1} F.\ Balibrea, J.C.\ Valverde, Bifurcations under non-degenerated conditions of 
higher degree and a new simple proof of the Hopf bifurcation theorem, J. Math. Anal. Appl., \textbf{237} (1999), 93-105.

\bibitem{bsv} H. Broer, C. Sim$\acute{o}$, R. Vitolo, The Hopf-saddle-node bifurcations for fixed points of 3D-diffeomorphisms: the Arnol'd resonance web, Bull.Belg.Math.Soc. Simon Stevin 15 (2008), no. 5, Dynamics in perturbations, 769-787. 

\bibitem{ca} J. Carr, \textit{Applications of Center Manifold Theory}, Springer-Verlag, New York, 1981

\bibitem{ch1} A. Chenciner, Bifurcations of elliptic fixed points. I. Invariant curves, Inst. Hautes \' Etudes Sci.Publ.Math. No. 61 (1985), 67-127

\bibitem{ch2}  A. Chenciner, Bifurcations of elliptic fixed points.II. Periodic Orbits and Invariant Cantor Sets, Invent. Math. 80 (1985), no.1 , 81-106

\bibitem{de} R.L. Devaney, \textit{Introduction to Chaotic Dynamical Systems}, The Benjamin/Cummings, New York, 1986

\bibitem{neveda} N.\ Elezovi\'c, V.\ \v Zupanovi\'c, D.\ \v Zubrini\'c, Box dimension of trajectories 
of some discrete dynamical systems, Chaos, Solitons \& Fractals Vol.\ 34,  2 (2007), 244-252. 

\bibitem{fa} K. Falconer, \textit{Fractal Geometry: Mathematical Foundations and Applications}, Chichester: John Wiley and Sons, USA, 1990

\bibitem{guho} J. Guckenheimer, P. Holmes, \textit{Nonlinear Oscillations, Dynamical systems, and Bifurcations of Vector Fields}, Springer-Verlag New York, USA, 1983

\bibitem{laho3} L. Horvat Dmitrovi\' c, Fractal Analysis of Bifurcations of Discrete Dynamical Systems and Applications to Continuous Systems, PhD thesis, University of Zagreb, 2011

\bibitem{laho1} L. Horvat Dmitrovi\' c, Box dimension and bifurcations of one-dimensional discrete dynamical systems, Discrete Contin. Dyn. Syst. 32 (2012), no. 4, 1287-1307. 

\bibitem{laho2} L. Horvat Dmitrovi\' c, Box dimension and hyperbolicity of fixed point/singularity of dynamical systems in $\mathbb{R}^{n}$, preprint 2012
          
\bibitem{hozu} L. Horvat Dmitrovi\' c, V. \v Zupanovi\' c, Box dimension of unit time map near nilpotent singularity of planar vector field, preprint 2012, arXiv:1205.5478[math.DS]
        
\bibitem{kuz} Yu.A. Kuznetsov, \textit{Elements of Applied Bifurcation Theory}, Springer-Verlag New York, USA, 1998

\bibitem{kuz2} Yu.A. Kuznetsov, H.G.E. Meijer, Remarks on Interacting Neimark-Sacker Bifurcations, J. Difference Equ. Appl. 12 (2006), no. 10, 1009-1035.
     
\bibitem{lapo} M.L. Lapidus, C. Pomerance, The Riemann Zeta-function and the one-dimensional Weyl-Berry conjecture for fractal drums, 
Proc. Lon. Math. Soc., \textbf{3} (1993), 66(1):41-69. 

\bibitem{marezu} P. Marde\v si\' c, M. Resman and V. \v Zupanovi\' c, Multiplicity of fixed points and growth of $\varepsilon$-neighbourhoods of orbits, J. Differ. Equ. 253 (2012), 2493-2514. 
 
\bibitem{ne} J. Neimark, On some cases of periodic motions depending on parameters, Dokl. Akad. Nauk SSSR 129, 736-739, 1959

\bibitem{pa} M. Pa\v si\' c, Minkowski-Bouligand dimension of solutions of the one-dimensional $p$-Laplacian,
        J. Differential Equations 190 (2003), 268-305
        
\bibitem{pazuzu} M. Pa\v si\' c, D. \v Zubrini\' c, V. \v Zupanovi\' c, Oscillatory and phase dimensions of solutions of some second-order differential equations, Bull. Sci. math. 133 (8) (2009), 859-874

\bibitem{pazuzu2} M. Pa\v si\' c, D. \v Zubrini\' c, V. \v Zupanovi\' c, Fractal properties of solutions of differential equations, chapter in a "Classification and Applications of Fractals", NOVA Science Publishers, 2011

\bibitem{pe} L. Perko, \textit{Differential Equations and Dynamical Systems}, Springer-Verlag New York, USA, 1996

\bibitem{sa} R. Sacker, A new approach to the perturbation theory of invariant surfaces, Comm. Pure Appl. Math. 18 (1965) 717-732

\bibitem{t} C. Tricot, \textit{Curves and Fractal dimension}, Springer-Verlag NY, 1995

\bibitem{wig} S. Wiggins, \textit{Introduction to Applied Nonlinear Systems and Chaos}, Springer-Verlag, New York, 1990

\bibitem{zu2} D. \v Zubrini\' c, Analysis of Minkowski contents of fractal sets and applications, Real Analysis Exchange, Vol 31(2), 2005/2006, 315-354

\bibitem{zuzu} D.\ \v Zubrini\'c, V.\ \v Zupanovi\'c, Fractal analysis of spiral trajectories of some
planar vector fields, Bull. Sci. math. 129/6 (2005), 457-485

\bibitem{zuzu3} D.\ \v Zubrini\'c, V.\ \v Zupanovi\'c,
Fractal analysis of spiral trajectories of some vector fields in ${\mathbb R}^3$, 
C.\ R.\ Acad.\ Sci.\ Paris, S\'erie I, Vol.\ 342, 12 (2006), 959-963

\bibitem{zuzu4} V.\ \v Zupanovi\' c, D.\ \v Zubrini\' c, \textit{Fractal dimension in dynamics}, in Encyclopedia of Math. Physics, J.-P. Fran\c coise, G.L. Naber, S.T. Tsou (Eds.), vol. 2 Elsevier, Oxford 2006.

\bibitem{belg} D.\ \v Zubrini\'c, V.\ \v Zupanovi\'c, Poincar\'e map in fractal analysis of 
spiral trajectories of planar vector fields, Bull.\ Belg.\ Math.\ Soc.\ Simon Stevin, 15(2008) 947-960

        
\end{thebibliography}
\end{document}